\documentclass{gtart_h}  

\def\ifplaintex{\expandafter\ifx\csname documentclass\endcsname\relax}

\ifplaintex 
\else
\fi

\expandafter\ifx\csname beginpicture\endcsname\relax
\expandafter\ifx\csname documentclass\endcsname\relax
\input pictex \else
\input prepictex \input pictex \input postpictex \fi\fi

\def\gt{{\mathsurround=0pt\it $\cal G\mskip-2mu$eometry \&\ 
$\cal T\!\!$opology}}        

\def\gtp{{\mathsurround=0pt\it $\cal G\mskip-2mu$eometry \&\ 
$\cal T\!\!$opology $\cal P\!$ublications}}  

\def\lognumber#1{\def\thelognumber{#1}}
\def\volumenumber#1{\def\thevolumenumber{#1}}
\def\papernumber#1{\def\thepapernumber{#1}}
\def\volumeyear#1{\def\thevolumeyear{#1}}

\def\pagenumbers#1#2{\def\startpage{#1}\def\finishpage{#2}}
\def\published#1{\def\publishdate{#1}}
\def\proposed#1{\def\theproposer{#1}}
\def\seconded#1{\def\theseconders{#1}}
\def\received#1{\def\receiveddate{#1}}
\def\revised#1{\def\reviseddate{#1}}
\def\accepted#1{\def\accepteddate{#1}}

\def\asciiaddress#1{\def\theasciiaddress{#1}}
\def\asciiemail#1{\def\theasciiemail{#1}}

\long\def\asciiabstract#1{\long\def\theasciiabstract{#1}}

\let\\\par\let\thelognumber\relax
\let\thevolumenumber\relax\let\thepapernumber\relax
\let\thevolumeyear\relax\let\thesamplenumber\relax\let\startpage\relax
\let\finishpage\relax\let\publishdate\relax\let\receiveddate\relax
\let\reviseddate\relax\let\accepteddate\relax\let\theasciititle\relax
\let\theasciiauthors\relax\let\theasciiaddress\relax
\let\theasciiabstract\relax
\let\theasciiemail\relax\let\theshortauthors\relax\let\theshorttitle\relax

\long\def\maketitlep{   

\count0=\startpage

\gt\hfill      
\beginpicture
\setcoordinatesystem units <0.33truein, 0.33truein> point at 2.2 0.9
\setplotsymbol ({$\cal G$})
\plotsymbolspacing=9truept
\circulararc 315 degrees from 0 1 center at 0 0
\setplotsymbol ({$\cal T$})
\circulararc 315 degrees from 1 -1 center at 1 0
\endpicture
\break
{\small\ifx\thesamplenumber\relax 
Volume \else Sample
\fi\thevolumenumber\ (\thevolumeyear)
\startpage--\finishpage\nl
Published: \publishdate}
\vglue 0.5truein plus 0.4fil minus 0.1truein

{\parskip=0pt\leftskip 0pt plus 1fil\def\\{\par\smallskip}{\ifplaintex\large
\else\Large\fi\bf\thetitle}\par\medskip}   

\vglue 0pt plus 0.1fil 

{\parskip=0pt\leftskip 0pt plus 1fil\def\\{\par}{\sc\theauthors}
\par\medskip}

\vglue 0pt plus 0.1fil 

{\small\parskip=0pt\let\newline\\
{\leftskip 0pt plus 1fil\def\\{\par}{\sl\theaddress}\par}
\expandafter\ifx\theemail\relax    
\relax\else\vglue 5pt plus 0.02fil minus 2pt\def\\{\stdspace{\rm 
and}\stdspace} 
\cl{Email:\stdspace\tt\theemail}\fi
\ifx\theurl\relax                  
\relax\else\vglue 5pt plus 0.02fil minus 2pt\def\\{\stdspace{\rm 
and}\stdspace}
\cl{URL:\stdspace\tt\theurl}\fi\par}

\vglue 7pt plus 0.3fil minus 3pt

{\bf Abstract}
\vglue 5pt plus 0.1fil minus 2pt

\theabstract

\vglue 7pt plus 0.3fil minus 3pt

{\bf AMS Classification numbers}\quad Primary:\quad \theprimaryclass

Secondary:\quad \thesecondaryclass

\vglue 5pt plus 0.3fil minus 2pt

{\bf Keywords:}\quad \thekeywords

\vglue 10pt plus 0.5fil minus 5pt

{\small  Proposed: \theproposer\hfill Received: \receiveddate\nl
Seconded: \theseconders\hfill 
\ifx\reviseddate\relax                         
Accepted: \accepteddate                        
\else
Revised: \reviseddate                          
\fi}
\eject
}       

\let\maketitlepage\maketitlep
\let\maketitle\maketitlepage

\font\phead=cmsl9 scaled 950
\font=cmsl9 scaled 1050
\font\pnum=cmbx10 scaled 913
\font\lnum=cmbx10 
\font\pfoot=cmsl9 scaled 950
\font=cmsl9 scaled 1050
\ifplaintex
\headline{\vbox to 0pt{\vskip -4.5mm\line{\small\phead\ifnum
\count0=\startpage ISSN 1364-0380 (on line)
1465-3060 (printed) \hfill {\pnum\folio}\else\ifodd\count0\def\\{ }%
\ifx\theshorttitle\relax\thetitle\else\theshorttitle\fi\hfill{\pnum\folio}
\else\def\\{ and }{\pnum\folio}\hfill\ifx\theshortauthors\relax\theauthors
\else\theshortauthors\fi\fi\fi}\vss}}
\footline{\vbox to 0pt{\vglue 0mm\line{\small\pfoot\ifnum\count0=\startpage
\copyright\ \gtp\hfill\else
\gt, Volume \thevolumenumber\ (\thevolumeyear)\hfill\fi}\vss
}}
\else
\makeatletter
\def\@oddhead{{\small\lhead\ifnum\count0=\startpage ISSN 1364-0380 (on line)
1465-3060 (printed) \hfill {\lnum\number\count0}\else\ifodd\count0
\def\\{ }\ifx\theshorttitle\relax \thetitle \else\theshorttitle\fi\hfill
{\lnum\number\count0}\else\def\\{ and }{\lnum\number\count0}
\hfill\ifx\theshortauthors\relax 
\theauthors\else\theshortauthors\fi\fi\fi}}\def\@evenhead{\@oddhead}
\def\@oddfoot{\small\lfoot\ifnum\count0=\startpage\copyright\ \gtp\hfill\else
\gt, Volume \thevolumenumber\ (\thevolumeyear)\hfill\fi}
\def\@evenfoot{\@oddfoot}
\makeatother
\fi

\newwrite\gtoutfile
\long\gdef\makeheadfile{  
{\def\\{, }\def\s{ }
\immediate\openout\gtoutfile head.xxx
\immediate\write\gtoutfile{Proxy-for: \ifx\theasciiauthors\relax
\theauthors\else\theasciiauthors\fi\s<\ifx\theasciiemail\relax\theemail\else\theasciiemail\fi>}
\immediate\write\gtoutfile{\noexpand\\}
\immediate\write\gtoutfile{Authors: \ifx\theasciiauthors\relax
\theauthors\else\theasciiauthors\fi}
{\def\\{ }\immediate\write\gtoutfile{Title: \ifx\theasciititle\relax
\thetitle\else\theasciititle\fi}}
\immediate\write\gtoutfile{Subj-class: GT or SG or MG etc}
\immediate\write\gtoutfile{MSC-class: \theprimaryclass\ifx\thesecondaryclass\relax\else, \thesecondaryclass\fi}
\immediate\write\gtoutfile{Journal-ref: Geom. Topol. \thevolumenumber
(\thevolumeyear) \startpage-\finishpage}
\immediate\write\gtoutfile{Comments: Published by Geometry and Topology at}
\immediate\write\gtoutfile{\s\s http://www.maths.warwick.ac.uk/gt/GTVol\thevolumenumber/paper\thepapernumber.abs.html}
\immediate\write\gtoutfile{\noexpand\\}
\immediate\write\gtoutfile{}
\ifx\theasciiabstract\relax
\immediate\write\gtoutfile{\theabstract}\else
\immediate\write\gtoutfile{\theasciiabstract}\fi
\immediate\write\gtoutfile{}
\immediate\write\gtoutfile{\noexpand\\}
\immediate\write\gtoutfile{}
\immediate\closeout\gtoutfile}}  

\def\maketitlepage{\maketitlep\makeheadfile}
\let\maketitle\maketitlepage

\lognumber{415}
\received{30 January 2004}
\volumenumber{9}\papernumber{42}\volumeyear{2005}
\pagenumbers{1881}{1913}   
\revised{20 September 2005}
\published{6 October 2005}
\accepted{20 September 2005}
\proposed{Steve Ferry}
\seconded{Ralph Cohen, Leonid Polterovich}

\usepackage{amssymb,amsmath,graphicx,amscd,psfrag}

\def\psfraga <#1,#2> #3#4{%
\psfrag {#3}{\smash{\rlap{\kern #1 \raise #2\hbox{#4}}}}}

\def\figref#1{\hyperlink{#1anchor}{Figure~\ref*{#1}}}
\def\anchor#1{\noindent\hypertarget{#1anchor}{\smash{$\phantom{99}$}}}

\numberwithin{equation}{section} 

\theoremstyle{plain} 
\newtheorem{thm}{Theorem}[section] 
 
\newtheorem{cor}[thm]{Corollary} 
\newtheorem{lem}[thm]{Lemma} 
\newtheorem{lemma}[thm]{Lemma} 
 
\newtheorem{prop}[thm]{Proposition} 
\newtheorem{theodef}[thm]{Theorem-Definition} 

\newcommand\theoref{Theorem~\ref} 
\newcommand\lemref{Lemma~\ref} 
\newcommand\propref{Proposition~\ref} 
\newcommand\corref{Corollary~\ref} 
\newcommand\defref{Definition~\ref}

\theoremstyle{definition} 
\newtheorem{defn}[thm]{Definition}

\newtheorem{rem}[thm]{Remark}{} 
\newtheorem{rems}[thm]{Remarks} 
\newtheorem{ex}[thm]{Example} 
\newtheorem{exs}[thm]{Examples} 
\newtheorem{con}[thm]{Construction}

\def\Im{\protect\operatorname{Im}}

\def\Int{\protect\operatorname{Int}}

\def\sign{\protect\operatorname{sign}} 
\def\pal{\protect\operatorname{al}}

\def\win{\protect\operatorname{win}}

\def\gl{\protect\operatorname{gl}} 
\def\lk{\protect\operatorname{lk}} 
\def\fr{\protect\operatorname{fr}}
\def\hq{\protect\operatorname{hq}}
 
\def\pt{\protect\operatorname{pt}}

\def\al{\protect\operatorname{alk}} 
 
\def\no{\protect\operatorname{\mathcal N_1}} 
\def\nd{\protect\operatorname{\mathcal N_2}} 
\def\ovno{\overline{\mathcal N}_1}
\def\ovnd{\overline{\mathcal N}_2}
\def\B{\protect\operatorname{\mathcal B}} 
\def\A{\mathbb A} 
\def\Indet{\protect\operatorname{Indet}} 
\def\NN{\protect\operatorname{\mathcal N}} 

\def\ga{\alpha} 
\def\gb{\beta} 
\def\gl{\lambda} 
\def\eps{\varepsilon} 
\def\gf{\varphi} 
\def\bor{\Omega} 

\def\Z{{\mathbb Z}} 
\def\Q{{\mathbb Q}} 
\def\R{{\mathbb R}} 
\def\P{{\mathbb P}} 
 
\def\NN{\mathcal N}

\def\1{\hbox{\rm\rlap {1}\hskip.03in{\rom I}}} 
\def\Bbbone{{\rm1\mathchoice{\kern-0.25em}{\kern-0.25em} 
{\kern-0.2em}{\kern-0.2em}I}} 

\def\p{\partial}

\def\bul{^{\scriptstyle{\scriptstyle{\bullet}}}} 
\def\wt{\widetilde} 
\def\wh{\widehat}

\def\ov{\overline} 

\begin{document} 

\title{Toward a general theory of linking invariants} 
\author{Vladimir V Chernov\\Yuli B Rudyak} 
\address{Department of Mathematics, 6188 Bradley Hall\\Dartmouth College,
  Hanover NH 03755-3551, USA\\{\rm and}\\Department of Mathematics, University of 
Florida\\358
  Little Hall, Gainesville, FL 32611-8105, USA} 
\asciiaddress{Department of Mathematics, 6188 Bradley Hall\\Dartmouth College,
  Hanover NH 03755-3551, USA\\and\\Department of Mathematics, University of 
Florida\\358
  Little Hall, Gainesville, FL 32611-8105, USA} 
\gtemail{\mailto{Vladimir.Chernov@dartmouth.edu}{\rm\qua 
and\qua}\mailto{rudyak@math.ufl.edu}}
\asciiemail{Vladimir.Chernov@dartmouth.edu, rudyak@math.ufl.edu}

\begin{abstract} 
Let $N_1, N_2, M$ be smooth manifolds with $\dim N_1 
+\dim N_2 +1= \dim M$ and let $\phi_i$, for $i=1,2$, be 
smooth mappings of $N_i$ to $M$ where
$\Im \phi _1\cap \Im \phi _2=\emptyset$. 
The classical linking number $\lk (\phi_1,\phi_2)$ is 
defined only when $\phi_{1*} [N_1]=\phi_{2*}[N_2]=0\in H_*(M)$. 

The affine linking invariant $\al$ is a generalization of $\lk$
to the case where $\phi_{1*} [N_1]$ or $\phi_{2*}[N_2]$ are not
zero-homologous.
In \cite{ChernovRudyak} we constructed the 
first examples of affine linking invariants 
of nonzero-homologous spheres in the spherical tangent 
bundle of a manifold, and showed that $\al$ is intimately related to 
the causality relation of wave fronts on manifolds.  In this paper we
develop the general theory.

The invariant $\al$ appears to be a universal Vassiliev--Goussarov
invariant of order $\leq 1$. In the case where $\phi_{1*}
[N_1]=\phi_{2*}[N_2]=0\in H_*(M)$, it is a splitting  of the classical
linking number into a collection of independent invariants.

To construct $\al$ we introduce a new pairing $\mu$ on the bordism groups
of spaces of mappings of $N_1$ and $N_2$ into $M$, not necessarily under
the restriction $\dim N_1+\dim N_2 +1= \dim M$.  For the
zero-dimensional bordism groups, $\mu$ can be related to the
Hatcher--Quinn invariant. In the case $N_1=N_2=S^1$, it is 
related to the  Chas--Sullivan string homology super Lie bracket, and
to the Goldman Lie bracket of free loops on surfaces.
\end{abstract}

\asciiabstract{Let N_1, N_2, M be smooth manifolds with dim N_1 + dim
N_2 +1 = dim M$ and let phi_i, for i=1,2, be smooth mappings of N_i to
M with Im phi_1 and Im phi_2 disjoint.  The classical linking number
lk(phi_1,phi_2) is defined only when phi_1*[N_1] = phi_2*[N_2] = 0 in
H_*(M).  The affine linking invariant alk is a generalization of lk to
the case where phi_1*[N_1] or phi_2*[N_2] are not zero-homologous.  In
arXiv:math.GT/0207219 we constructed the first examples of affine
linking invariants of nonzero-homologous spheres in the spherical
tangent bundle of a manifold, and showed that alk is intimately
related to the causality relation of wave fronts on manifolds.  In
this paper we develop the general theory.  The invariant alk appears
to be a universal Vassiliev-Goussarov invariant of order < 2. In the
case where phi_1*[N_1] and phi_2*[N_2] are 0 in homology it is a
splitting of the classical linking number into a collection of
independent invariants.  To construct alk we introduce a new pairing
mu on the bordism groups of spaces of mappings of N_1 and N_2 into M,
not necessarily under the restriction dim N_1 + dim N_2 +1 = dim M.
For the zero-dimensional bordism groups, mu can be related to the
Hatcher-Quinn invariant.  In the case N_1=N_2=S^1, it is related to
the Chas-Sullivan string homology super Lie bracket, and to the
Goldman Lie bracket of free loops on surfaces.}


{\small\maketitlepage}

\section{Introduction} 

In this paper the word ``smooth'' means $C^{\infty}$. Throughout this paper $M$ 
is a smooth connected oriented 
manifold (not necessarily compact), and $N_1$, $N_2$ are 
smooth oriented closed manifolds. 
The dimensions of $M, N_1, N_2$ are denoted by $m, n_1, n_2 
$, respectively, and the one-point space is denoted by $\pt$. 

Let $\mathcal N_i$, for $i=1,2$, be a path-connected component of
the space of all smooth mappings of $N_i$ to $M$. (Thus the mappings
in $\mathcal N_i, i=1,2,$ are not assumed to be immersions.)  Let
$\B=\B_{\no,\nd}$ be the space of quadruples $(\phi_1, \phi_2, \rho_1,
\rho_2)$ where $\phi_i\co  N_i \to M, i=1,2,$ belong to $\mathcal N_i$ and
$\rho_i\co  \pt \to N_i$ are such that $\phi_1\rho_1=\phi_2\rho_2$. Clearly,
$\B$ can be regarded as a subset of $\no\times \nd \times N_1\times N_2$,
and we equip $\B$ with the subspace topology.

The classical linking number $\lk$ is a $\Z$--valued invariant
of a pair $(\phi_1, \phi_2)\in \no \times\nd$ with $n_1+n_2+1=m$
(and with $\phi_1(N_1)\cap \phi_2(N_2)=\emptyset$). The invariant
$\lk(\phi_1, \phi_2)$ is defined only if $\phi_{1*}([N_1]),
\phi_{2*}([N_2])=0\in H_*(M)$ (or if $\phi_{1*}([N_1]), \phi_{2*}([N_2])$
are torsion classes, in which case $\lk$ takes values in $\Q$ or
$\Q/\Z$). U~Kaiser~\cite{Kaiserbook} generalized linking numbers to
the case of arbitrary submanifolds of the linking dimension that are
homologous into a boundary or into an end of the ambient manifold.
(For $M$ being the solid torus the similar approach to defining linking
numbers was previously used by S~Tabachnikov~\cite{Tabachnikov}.)

The main goal of the paper is to construct a version of the linking
invariant $\lk$ for pairs $(\phi_1, \phi_2)\in \no \times\nd$ with
$n_1+n_2+1=m$ and $\phi_1(N_1)\cap \phi_2(N_2)=\emptyset$ without any
restrictions on the homology classes $\phi_{1*}([N_1]), \phi_{2*}([N_2])$.

In greater detail, let $\Sigma \in \no \times \nd$ be the subset
of pairs $(\phi_1, \phi_2)$ with $\phi_1(N_1)\cap \phi_2(N_2) \ne
\emptyset$. Fixing a pair $*\in \no \times \nd \setminus \Sigma$,
we define an invariant 
$$\al\co  \no \times \nd \setminus \Sigma\to H_0(\B)/\Indet$$ 
which is an invariant under link homotopy of pairs $(\phi_1, \phi_2)$; here 
$\Indet$ 
is a certain indeterminacy subgroup. We call $\al$ the {\it affine} linking 
invariant, since the change of the base point $*$ leads to changing of
$\al$ 
by an additive constant. 

It turns out that the augmentation $\B \to \pt$ reduces our invariants to 
the classical ones (ie the linking numbers with values in 
$H_0(\pt)=\Z$) provided that the last ones are defined. In other words, here 
we have a splitting of the classical linking invariant. 

Our constructions can be easily modified to yield affine linking type
invariants under the {\it singular concordance} relation. In the case
of 1--links in 3--manifolds this will give us the invariants constructed
by Schneiderman in~\cite{Schneiderman}.

Our construction can be rather easily modified to give an invariant of a pair 
$(\phi_1, \phi_2) \in \no \times \nd$ with 
$\phi_1(N_1)\cap 
\phi_2(N_2)=\emptyset$ 
and without any restrictions on the dimensions $n_1, n_2, m$ of $N_1, N_2, M$. 
In this case we also get an invariant of a link $(\phi_1, \phi_2)$ considered up 
to the link homotopy. 
The invariant takes values in 
$\bor_{n_1+n_2+1-m}(\B)/\bigl (\Im \mu_{1,0}+\Im \mu_{0,1}\bigr)$, where 
$\mu_{i,j}\co \bor_i(\no)\otimes \bor_j(\nd)\to \bor_{i+j+n_1+n_2-m}(\B)$ 
is the pairing defined in Theorem~\ref{pairing}. We plan to study this 
more general invariant in detail in another paper.

Most of our results are based on considering of a helpful pairing 
$$ 
\mu_{i,j}\co  \bor _i(\no )\otimes \bor_j (\nd ) \to \bor_{i+j+n_1+n_2-m} (\B ) 
$$ 
where $\bor_*(X)$ is the group of oriented bordisms of $X$. This pairing has 
many remarkable properties. For example:

(1)\qua The pairings $\mu_{1,0}, \mu_{0,1}$ enable us to describe the indeterminacy 
for the 
invariant $\al$. (Note that $H_0(X)=\bor_0(X)$ for all $X$.) 

(2)\qua The pairing $\mu_{0,0}$ tells us (in many cases) whether two 
$C^{\infty}$ maps $f_1\co  N_1 \to M$ and $f_2\co  N_2\to M$ can be deformed 
to maps with disjoint images, see Section~\ref{disjointmappings}.  The case when 
two {\it immersions} $f_1$ and 
$f_2$ can be {\it regularly} homotoped to maps with disjoint images was 
considered by Hatcher and Quinn~\cite{HQ}. Concerning relations between 
$\mu_{0,0}$ and  the Hatcher--Quinn invariant, see  subsection~\ref{hq}. A 
conicidence 
problem for the case $N_1=N_2$ was considered by 
Koschorke~\cite{Koschorke1} via the approach of Hatcher--Quinn invariants.

(3)\qua If $N_1=N_2=S^n$ and $M$ is a $2n$--manifold then $\mu$ leads to a 
generalization 
of  the Goldman bracket \cite{Goldman} of free loops on $2$--surfaces, see 
 subsection \ref{goldman}. 

(4)\qua In case of $N_1=N_2=S^1$ the pairing  $\mu$ leads to a (graded) Lie algebra 
structure on $\bor_*(\ov \NN)$ where $\ov \NN=\ovno=\ovnd$ is the union of 
all the connected components of the space of mappings $S^1\to M$. This Lie 
algebra structure 
is 
related to
the string homology Lie bracket introduced by 
Chas and Sullivan~\cite {ChasSullivan},~\cite{ChasSullivan2}, cf~also the work 
of 
A~Voronov~\cite{Voronov}.
We are not able to discuss this algebra in detail here but intend to do it in 
the coming development of our work~\cite{ChernovRudyakStrings}. 

(5)\qua In fact the mapping $\mu$ extends to a Lie bracket on the nonoriented bordism 
groups of
mappings into $M$ of garlands glued out of arbitrary manifolds.
It also seems that for the appropriately chosen signs $\mu$ extends to a (super) 
Lie bracket even for oriented bordism groups, but we are still computing the 
appropriate signs in the graded Jacobi identity, see~\cite[Theorem 
3.1]{ChernovRudyakStrings}.

The paper is organized as follows. In Section \ref{secbor} we introduce the 
pairing ${\mu}$. In Sections \ref{affinesection} and \ref{proofmainthm} we 
define affine linking invariants of pairs $(\phi_1, \phi_2)\in \no \times\nd$ 
(with 
$f_1(N_1)\cap f_2(N_2)=\emptyset$) as elements of the group $\bor_0(\B)$ modulo 
certain indeterminacy; the last one is described in terms of the pairing $\mu$. 
In Section \ref{secrelation} we prove that the augmentation $\bor_0(\B)=H_0(\B) 
\to H_0(\pt)=\Z$ induced by the mapping $\B\to \pt$ reduces 
our invariants to the classical linking invariant, when the last one is 
well-defined. In Section \ref{Indet=0} 
we give conditions that guarantee the vanishing of the indeterminacy. In 
Section \ref{B} we give an explicit description of $\pi_0(\B)$ and $\bor_0(\B)$. 
In Section \ref{disjointmappings} we show that, in many cases, the pair $\bar 
f_1\co  N_1 \to M$, $\bar f_2\co  N_2 \to M$ of maps, with $n_1+n_2=m$ is homotopic to 
a pair $(f_1,f_2)$ with disjoint images, $f_1(N_1) \cap f_2(N_2)=\emptyset$, if 
and only if the pairing $\mu$ takes the zero value on $(\bar f_1, \bar f_2)$. If 
$\bar f_1$ and  $\bar f_2$ are homotopic to immersions, then the results of 
Section~\ref{disjointmappings} follow immediately from Theorem~$2.2$ of 
Hatcher--Quinn~\cite{HQ}. 

\paragraph{Acknowledgements} The first author was partially supported
by the Walter and Constance Burke Research Initiation Award. The second
author was partially supported by NSF, grant 0406311, and by MCyT,
projects BFM 2002-00788 and MCyT BFM2003-02068/MATE, Spain; his visit to
Dartmouth College was supported by the funds donated by Edward Shapiro
to the Mathematics Department of Dartmouth College.

We are also grateful to the anonymous referee for the very valuable comments and 
suggestions.

\section[The pairing in bordism]{The pairing $\mu_{i,j}\co  \bor _i(\no )\otimes \bor_j 
(\nd ) 
\to \bor_{i+j+n_1+n_2-m} (\B )$}\label{secbor} 
In this section we do not assume that $\dim N_1+\dim N_2+1=\dim M$.

Given a space $X$, we denote by $\bor_n(X)$ the $n$--dimensional 
oriented bordism group of $X$. Recall that $\bor_n(X)$ is 
the set of the 
equivalence classes of (continuous) maps $f\co  V^n \to X$ 
where $V$ is a 
closed oriented manifold. Here two maps $f_1\co  V_1 \to X$ 
and $f_2\co  V_2 
\to X$ are equivalent if there exists a map $g\co  W^{n+1}\to 
X,$ where $W$ 
is a compact oriented manifold whose oriented boundary $\p 
W$ is 
diffeomorphic to $V_1\sqcup (-V_2)$ and $g|_{\p W}=f_1\sqcup f_2$. 
Furthermore, the 
operation of disjoint union converts $\bor_n(X)$ into an 
abelian group. 
See \cite{Rudyak, Stong, Switzer} for details. 

Let $[V]\in H_n(V)$ be the fundamental class of a closed 
oriented 
$n$--dimensional manifold $V$. Every map $f\co  V \to X$ gives 
us an 
element $f_*[V]\in H_n(X)$, and the correspondence $(V,f) 
\mapsto 
f_*[V]$ yields the Steenrod--Thom homomorphism 
$$ 
\tau\co  \bor_n(X) \to H_n(X). 
$$ 
It turns out \cite{Thom} that this homomorphism is an isomorphism 
for $n\le 3$ 
and an epimorphism for $n\le 6$, see \cite 
{Rudyak, Stong} 
for modern proofs. 

Let $\alpha_1\co F_1\to \no$ be a mapping representing 
$[\alpha_1]\in \bor _i(\no)$ and let $\alpha_2\co F_2\to 
\nd$ be a mapping representing 
$[\alpha_2]\in \bor_{j}(\nd)$. 
Let $\wt \alpha_i\co F_i\times N_i\to M$, $i=1,2$, be such 
that $\wt \alpha_{i}(f,n)=(\alpha_{i}(f))(n)$. 
Let $v_1\in F_1\times N_1$ and $v_2\in F_2\times N_2$ be 
such that $\wt \alpha_1 (v_1)=\wt \alpha_2 (v_2)$. We say 
that 
$\wt \alpha_1$ and $\wt \alpha_2$ are {\it transverse} at $(v_1, 
v_2)$ if 
$d \wt \alpha_2\bigl( T_{v_2}(F_2\times N_2)\bigr)$ and 
$d\wt \alpha_1 \bigl (T_{v_1}(F_1\times N_1)\bigr )$ 
generate $T_{\wt \alpha_1(v_1)}M=T_{\wt \alpha_2(v_2)}M$. 
Following standard arguments we can assume that 
$\wt \alpha_1$ and $\wt \alpha_2$ are transverse, 
ie they are transverse at all $(v_1, v_2)$ such that 
$\wt \alpha_1 (v_1)=\wt \alpha_2 (v_2)$. 

Consider the pullback diagram 
\begin{equation}\label{bordism} 
\begin{CD} 
V@>j_1>> F_1\times N_1\\ 
@VVj_2V @VV\wt \alpha _1 V\\ 
F_2\times N_2 @>\wt \alpha _2 >> M\\ 
\end{CD}
\end{equation} 
of the maps $\wt \alpha_i$, for $i=1,2$. 

\begin{lem}\label{pullback} 
If $\wt \alpha_1,\wt \alpha_2$ are transverse, then 
$V$ is a closed oriented
$(i{+}j{+}$ $n_1{+}n_2{-}m)$--dimensional manifold.\qed
\end{lem} 

Let $p_1\co  F_i\times N_i \to F_i$ and $p_2\co  F_i\times N_i 
\to N_i$, for $i=1,2$, be the 
obvious projections. Consider the mapping 
$$ 
\mu (\wt \alpha_1, \wt \alpha_2)\co V\to 
\mathcal B,\quad v\mapsto \bigl(\phi_1^v, 
\phi_2^v,\rho_1^v, \rho_2^v) 
$$ 
where 
$\phi_i^v(n)= \wt \alpha_i(p_1(j_i(v)),n)$ for 
$n\in N_i$, and $\rho_i^v(\pt)=p_2(j_i(v))$. 

\begin{thm}\label{pairing} 
The above construction yields a well-defined pairing 
\begin{eqnarray*} 
\mu=\mu_{ij}\co  \bor _i(\no )\otimes \bor_j(\nd)&\to&
  \bor_{i+j+n_1+n_2-m} (\B),\\ 
\mu \left([\alpha_1], [\alpha_2]\right)&=&[V, \mu (\wt \alpha_1,
  \wt \alpha_2)].
\end{eqnarray*} 
\end{thm} 

\begin{proof} The routine argument shows that the bordism class $[V, \mu (\wt 
\alpha_1, 
\wt \alpha_2)]$ in 
$\bor _{i+j+n_1+n_2-m} (\B )$ depends only on $[\alpha_1]\in \bor _i (\no)$, 
$[\alpha_2]\in \bor _j(\nd)$. The bilinearity of $\mu$ is obvious. 
\end{proof}

\subsection{Relation to the Goldman bracket}\label{goldman} 
Let $\wt \NN_1$ (respectively, $\wt \NN_2$) be the topological space
that is the union of all the
connected components of the space of mappings of $N_1$ (respectively,
$N_2$) into $M$. 
Let $\wt \B$ be the topological space that is the union of $\B_{\NN_1, \NN_2}$ 
over all the connected components $\NN_1, \NN_2$.

Let $\ov \NN_i$ be the closure of $\wt \NN_i$ in the space of 
all continuous maps. Similarly, define $\ov \B$ to be the 
closure of $\wt \B$. 

Clearly, the pairing $\mu$ can be extended to 
$$\ov \mu\co  \bor _i(\ovno )\otimes \bor_j(\ovnd )\to 
\bor_{i+j+n_1+n_2-m}(\ov \B).$$
Let $M$ be a $2n$--dimensional manifold and let $N_1=N_2=S^n$. Given two
points $x,y\in S^n$ (not necessarily distinct), choose an orientation
preserving diffeomorphism $u_{x,y}\co  S^n\to S^n$ that maps $x$ to $y$.

For $M^{2n}$ and $N_1=N_2=S^n$ put 
$\ov \NN=\ovno=\ovnd$ and consider the pairing 
\begin{equation}\label{ov} 
\ov \mu\co  \bor_0(\ov\NN) \otimes \bor_0(\ov\NN) \to \bor_0 
(\ov\B). 
\end{equation} 
Since $\bor_0=H_0$, the pairing \eqref{ov} yields a pairing 
$$ 
\CD 
H_0(\ov\NN) \otimes H_0(\ov\NN) @>\ov\mu >>H_0 
(\ov\B). 
\endCD 
$$ 
Every point $b=(\phi_1, \phi_2, \rho_1, \rho_2)\in \ov \B$ gives us 
a map $h_b\co  S^n\vee S^n \to M$ as follows. We regard the sphere $S^n$ as a
pointed space with the base point $*$. Clearly, both maps 
$\phi_i u_{*,\rho_i(\pt)}\co  S^n\to M$, for $i=1,2$, map the point $*$ to 
the same point of $M$ and therefore yield the map $S^n\vee S^n \to M$.

The pairing $\mu\co  H_0(\ov\NN)\otimes H_0(\ov\NN)\to H_0(\ov\B)$ 
has the following interpretation. 
Put $\widehat \pi_n$ to be the orbits of $\pi_n(M)$ under the natural action 
of $\pi_1(M)$. Then $H_0(\ov \NN)$ is 
naturally identified with $\Z\widehat \pi_n$, a free $\Z$--module over 
$\widehat \pi_n$. Observe that any two modules $\Z\widehat \pi_n(M,p)$ 
and $\Z\widehat \pi_n(M,q)$ are 
canonically isomorphic

Let $s_i\co S^n\to M^{2n}$, for $i=1,2$, be 
two smooth generic mappings transverse to each other and realizing 
$[s_i]\in\widehat \pi_n$. Here genericity means that each intersection 
point of $s_1$ and
$s_2$ is {\em not} a self-intersection point of $s_1$ or $s_2$.
Put $P=\Im (s_1)\cap \Im (s_2)$.
Then 
\begin{equation}\label{muformula}
\ov \mu(1[s_1]\otimes 1[s_2])=\sum_{p\in P} 
\sign (p) [s_1,s_2, p], 
\end{equation}
where $\sign p$ is the natural orientation of $p$ coming from the 
intersection pairing. Here $[s_1,s_2, p]\in \pi_0(\ov \B)$ is the element 
that  maps the first sphere as $s_1$, the second sphere as $s_2$, 
and maps $\pt$ to the preimages of $p$ on the two spheres.

Now, using the coproduct $S^n\to S^n 
\vee S^n$, we conclude that every $b\in \ov \B$ gives us a map 
$$
S^n \longrightarrow S^n \vee S^n \stackrel{h_b}{\longrightarrow} M
$$
So, we get a map $\varphi\co  \ov \B \to \ov \NN$. Notice that this map is 
{\em not} continuous, but it induces a well-defined map 
$\pi_0(\ov \B) \to \pi_0(\ov \NN)$ since each $u_{x,y}$ is homotopic 
to the identity. Furthermore, since $\bor_0=H_0$, the pairing \eqref{ov} yields 
a pairing 
$$ 
\alpha\co  H_0(\ov\NN) \otimes H_0(\ov\NN)
\stackrel{\ov\mu}{\longrightarrow} H_0 
(\ov\B) \stackrel{\gf_*}{\longrightarrow} H_0(\ov\NN). 
$$ 
Now, because of the equality \eqref{muformula}, we conclude that
\begin{equation}
\label{spheresformula}
\alpha (1[s_1]\otimes 1[s_2])=\Bigl (\sum_{p\in P} 
\sign (p) [s_1s_2\in \pi_n(M, p)]\Bigr)\in \Z\widehat \pi_n. 
\end{equation}
Here the element $[s_1s_2\in \pi_n(M, p)]$ is the class in 
$\widehat \pi_n$ that is 
the product of the elements of 
$\pi_n(M,p)$ realized by $s_1$ and $s_2$.
More accurately, we have to compose $s_i$ 
with an automorphism $u_{*,x_i}$ of $S^n$ that maps the base point $*$
to the preimage $x_i=s^{-1}_i(p)$ of $p$. 

For $n=1$ the action of $\pi_1(M)$ on $\pi_n(M)=\pi_1(M)$ is given by 
the conjugation. 
So $\Z\widehat \pi_n=\Z\widehat \pi_1$ is 
a free $\Z$--module generated by the homotopy classes of free loops on $M$. 
Because of the explanations above $[s_1s_2\in \pi_n(M, p)]$ is a
well-defined element of $\Z\widehat \pi_1$.

Now formula 
$\eqref{spheresformula}$ is identical to the definition of
Goldman's Lie bracket on 
$\Z \widehat \pi_1(M^2)$, see Goldman~\cite[page 267]{Goldman}. 
Since both $\alpha$ and 
Goldman's Lie bracket are bilinear and coincide on the generators of 
$\Z \widehat \pi_1$, they are equal.
Thus for $N_1=N_2=S^1$ and $M^2$ our pairing $\mu$ generalizes 
Goldman's Lie bracket discovered 
in stages by Goldman~\cite{Goldman} and Turaev~\cite{Turaev}. 

\subsection{Relation to the Hatcher--Quinn invariant}\label{hq} 
Consider two maps $f\co  N_1 \to M$ and $g\co  N_2 \to M$ in $\no$ and $\nd$, 
respectively. 
Hatcher and Quinn~\cite{HQ} considered the homotopy pullback diagram 
$$ 
\CD E(f,g) @>f_E>>N_1\\ 
@Vg_E VV @VVfV\\ 
N_2 @>g>> M 
\endCD 
$$ 
For $i=1,2$, let $\nu_i$ be the stable normal vector bundle over $N_i$ and let 
$\xi$ be the stable vector bundle $f_E^*\nu_1\oplus g_E^*\nu_2\oplus 
f_E^*f^*\tau_M$. 
Let $E$ denote $E(f,g)$ and let $\Omega_k^{\fr}(E,\xi)$ denote 
the bordism group based on singular manifolds 
$h\co V^k \to E$ together with a stable bundle isomorphism of 
the normal bundle of $V$ with $h^*\xi$. 
The  $\Omega_k^{\fr}(E,\xi)$ groups are the homotopy groups 
of the Thom spectrum $T\xi$. 

Given two transverse maps $f_1\co  N_1 \to M$ 
and $f_2\co  N_2\to M $ homotopic to $f$ and $g$ respectively, consider the 
pullback $V$ 
(rather than the homotopy pullback) 
of $f_1$ and $f_2$. 
Consider the obvious maps $v_i\co V\to N_i$ and construct a map
$h\co  V \to E$ 
with $f_E\circ h= v_1$ and $g_E\circ h = v_2$. In fact, this $h$
determines and 
is determined by homotopies from $f_1$ to $f$ and from $f_2$ to $g$. 
Then $(V,h)$ yields an element of $\Omega_{n_1+n_2-m}^{\fr}(E,\xi)$,
and this element does not depend on the choice of the above described
homotopies. 
So, we have the Hatcher--Quinn map 
$$ 
\hq\co  \pi_0(\no) \times \pi_0(\nd)\longrightarrow 
\Omega^{\fr}_{n_1+n_2-m}(E,\xi). 
$$
Now, assume that $n_1+n_2=m.$ Take a point of $\B$ and represent it by
a commutative diagram 
$$ 
\CD 
\pt @>>>N_1\\
@VVV @VVf_1V\\
N_2 @>f_2>> M 
\endCD 
$$  
with $f_1, f_2$ transverse and $f_1\simeq f,\,f_2\simeq g$. Clearly,
this diagram gives us an element of the group 
$\Omega_{0}^{\fr}(E,\xi)$, 
and in fact we have a map
$\gf\co  \bor_{0}(\B)\to \Omega_{0}^{\fr}(E,\xi)$. It 
is easy to see that the diagram 
$$ 
\CD \pi_0(\no) \times \pi_0(\nd) @>\hq>>\Omega_{0}^{\fr}(E,\xi)\\ 
@VVV @AA\gf A\\ 
\bor_0(\no)\otimes\bor_0(\nd) @>\mu_{0,0}>> \bor_{0}(\B) 
\endCD 
$$ 
commutes.  So in the case $n_1+n_2=m$ the pairing $\mu_{0,0}$ can be
regarded as a version of the Hatcher--Quinn invariant.

\section{Affine linking invariants}\label{affinesection} 

From here and till the end of the paper we assume that 
$n_1+n_2+1=m$ (unless the opposite is explicitly stated). 

Put $\Sigma$ to be the discriminant in $\no \times 
\nd$, 
ie the subspace that consists of pairs $(f_1, f_2)$ such 
that 
there 
exist $y_1\in N_1, y_2\in N_2$ with $f_1(y_1)=f_2 (y_2)$. 
(We do not include into $\Sigma$ the maps that are singular 
in the 
common sense but do not involve double points between $f_1 
(N_1)$ and $f_2(N_2)$.) 

\begin{defn}\label{sigma0} 
We define $\Sigma_0$ to be the subset (stratum) of $\Sigma$ 
consisting 
of all 
the pairs $(f_1, f_2)$ such that there exists precisely one pair ($y_1, y_2)$ 
of points $y_1\in N_1, y_2\in N_2$ such that
\begin{enumerate}
\item[(a)] $f_1(y_1)=f_2(y_2)$;
\item[(b)] $y_i$ is a regular point of $f_i$, for $i=1,2$; 
\item[(c)] $(df_1)(T_{y_1}N_1)\cap (df_2)(T_{y_2}N_2)=0$. 
\end{enumerate}
\end{defn} 

\begin{con}\label{vector} 
Let $\gamma\co (a,b)\rightarrow \no \times \nd$ be a path 
which intersects $\Sigma_0$ in a point $\gamma(t_0)$. We also 
assume that 
$$\gamma(t_0-\delta, t_0+\delta)\cap \Sigma =\gamma(t_0)$$ 
for $\delta$ small enough. We construct a vector ${\mathbf 
v}={\mathbf v}(\gamma,t_0, \delta)$ as follows. 
We regard $\gamma(t_0)$ as a pair $(f_1,f_2)\in \mathcal 
\no 
\times \nd$ 
and consider the points $y_1,y_2$ as in Definition~\ref{sigma0}. Set 
$z=f_1(y_1)=f_2(y_2)$. 
Choose a small $\delta>0$ and regard $\gamma(t_0+\delta)$ 
as a pair $(g_1,g_2)\in \no \times \nd$. Set 
$z_i=g_i(y_i), i=1,2$. 
Take a chart for $M$ that contains $z$ and $z_i,\, i=1,2,$ 
and set 
$${\mathbf v}(\gamma,t_0,\delta) := 
  \overrightarrow{zz_1}-\overrightarrow{zz_2}\in T_z M.$$ 
\end{con} 

\begin{defn}\label{transversal} 
Let $\gamma\co (a,b)\rightarrow \no \times \nd$ be a path 
as in Construction~\ref{vector}. We say that $\gamma$ {\it intersects 
$\Sigma_0$ 
transversally for $t=t_0$} if there exists $\delta_0>0$ 
such that 
\[ 
{\mathbf v}(\gamma, t_0, \delta) \notin (df_1)(T_{y_1}N_1) 
\oplus 
(df_2)(T_{y_2}N_2) \subset T_zM 
\] 
for all $\delta \in (0, \delta_0)$. 
\end{defn} 

It is easy to see that the concept of transverse 
intersection 
does not depend on the choice of the chart. 

\begin{defn}\label{genericpath} 
A path $\gamma\co (a,b)\rightarrow \no\times \nd, -\infty 
\le a<b\le \infty$ is said to be {\it generic} if 
\begin{enumerate} 
\item[(a)] $\gamma(a,b)\cap \Sigma =\gamma(a,b)\cap \Sigma_0$; 
\item[(b)] the set $J=\{t|\gamma(t)\cap \Sigma_0\ne \emptyset\} \subset (a,b)$ 
  is an isolated subset of $\R$; 
\item[(c)] the path $\gamma$ intersects $\Sigma_0$ transversally for
  all $t\in J$.
\end{enumerate} 
\end{defn} 

As one can expect, every path can be turned into a generic 
one by a small deformation. We leave the proof to the 
reader. 

\begin{defn}\label{signalk} 
Let $\gamma$ be a path in $\no \times \nd$ that intersects 
$\Sigma$ transversally in one point $\gamma(t_0)\in \Sigma_0$. We 
associate a sign $\sigma(\gamma, t_0)$ to such a crossing as follows. 

We regard $\gamma(t_0)$ as a pair 
$(f_1,f_2)\in \no \times \nd$ and consider the points 
$y_1\in N_1,y_2\in N_2$ such that $f_1(y_1)=f_2(y_2)$. 
Set $z=f_1(y_1)=f_2(y_2)$. 
Let $\mathfrak r_1$ and $\mathfrak r_2$ be frames which are 
tangent to $N_1$ and $N_2$ at $y_1$ and $y_2$, respectively, 
and both are assumed to be positive. Consider the frame 
$$\{ df_1(\mathfrak r_1), {\mathbf v}, df_2(\mathfrak r_2) \}$$ 
at $z\in M$, where ${\mathbf v}$ is a vector described in 
Construction~\ref{vector}. We put $\sigma(\gamma, t_0)=1$ if this frame 
gives us the orientation of $M$, otherwise we put $\sigma(\gamma,t_0)=-1$. 
Because of the transversality and condition (c) from Definition~\ref{sigma0}, 
the family $ \{ df_1(\mathfrak r_1), {\mathbf v}, df_2 
(\mathfrak r_2) \} $ is indeed a frame. 
Note also that the vector ${\mathbf v}$ is not well-defined, 
but the above defined sign $\sigma$ is. 
\end{defn} 

Clearly if we traverse the path $\gamma$ in the opposite 
direction then the sign of the crossing changes.

For every space $X$, the group $\bor_0(X)=H_0(X)$ is the 
free abelian 
group with the base $\pi_0(X)$. So, every element of 
$\bor_0 
(X)$ can be 
represented as a finite linear combination $\sum \gl_kP_k$ 
with $\gl_k\in 
\Z$ and $P_k\in X$, and every such linear combination gives 
us an element 
of $\bor_0(X)$. Below in Section~\ref{B} we give examples 
of situations 
where $\pi_0(\B)$ is an infinite set, in spite of the fact 
that $\no$ and $\nd$ 
are path connected.

\begin{defn}\label{pointbordism} 
Let $\gamma$ be a path in $\no \times \nd$ that intersects 
$\Sigma$ transversally in one point $\gamma(t_0)\in 
\Sigma_0 
$. We define $[\gamma(t_0)]\in \bor_0(\B )$ as $\sigma(t_0) 
\gamma(t_0)$. 

Clearly, 
\begin{equation}\label{eps} 
\eps([\gamma(t_0])=\sigma(\gamma, t_0), 
\end{equation} 
where $\eps\co  \bor_0(\B) \to \Z=\bor _0(\pt)$ is the 
homomorphism induced 
by the map 
$\B \to \pt$. 
\end{defn} 

\begin{defn}\label{a} 
Since $\bor_i(X)=H_i(X)$ for $i=0,1$, by the
K\"unneth formula we have the canonical isomorphism 
$$ 
\bor_1(X \times Y)= \bor_1(X)\otimes \bor_0(Y) \oplus \bor_0 
(X)\otimes 
\bor_1(Y). 
$$ 
Now, the pairings 
$$ 
\mu_{10}\co \bor_1(\no)\otimes \bor_0(\nd) \to \bor_0(\B ) \text 
{ and } 
\mu_{01}\co \bor_0(\no)\otimes \bor_1(\nd) \to \bor_0(\B ) 
$$ 
yield the homomorphism 
$$ 
\gl\co  \bor_1(\no \times \nd) \to \bor_0(\B ). 
$$ 
We define the {\it indeterminacy subgroup} $\Indet\subset 
\bor_0(\B)$ to 
be the image of the homomorphism $\gl$. We define the 
abelian group $\A=\A(\no ,\nd)$ to be the 
quotient group of 
$\bor_0(\B)/\Indet$. Let $q=q_{\no, \nd }\co \bor_0(\B)\to \A$ 
be the 
corresponding quotient homomorphism. 
\end{defn} 

\begin{exs}
(1)\qua Consider the case of $M^3$ being a lens space with a fundamental group 
$\Z_p$, 
$N_1=S^2, N_2=\pt$. In this case there is only one homotopy class of mappings 
$S^2\to M$ and of $\pt \to M$. It is easy to see that $\bor_0(\B)=\Z$.

Take $f\co S^2\to M$ and $g\co \pt\times S^1=S^1\to M$. Then $\mu_{0,1}(f,g)$ equals 
to the intersection index between $f_*[S^2]$ and $g_*[S^1]\in H_*(M)$. Since 
$H_1(M)=\Z_p$ we get that $\Im \mu_{0,1}=0\in \Z=\bor_0(\B)$.

Take $r\co S^2\times S^1\to M$ and $s\co \pt \to M$. Then $\mu_{1,0}(r,s)\in \Z$ 
equals to the degree of the mapping $r\co S^2\times S^1\to M$. Since $\pi_2(M)=0$, 
the elementary obstruction theory shows that for a given image of 
$r_*(\cdot\times S^1)\in \pi_1(M)$ all the homotopy classes of mappings 
$S^2\times S^1\to M$ are classified by $\pi_3(M)$. Thus $\Im \mu_{1,0}$ 
coincides with the possible degrees of mappings $S^3\to M$. Since all such 
mappings pass through the covering map $S^3\to M$ that has degree $p$, we get 
$\Im \mu_{1,0}=p\Z\subset \Z$. 
Thus $\Indet=p\Z\subset \Z=\bor_0(\B)$.

(2)\qua A harder example comes from $M=F_g\times S^1$, where $F_g$
is an oriented surface of genus 
$g>1$. Let $N_1=N_2=S^1$ and let $\alpha, \beta\co S^1\to M$ be linked embedded 
circles that project to two simple curves on $F_g$ with one intersection point. 
Let $\no, \nd$ be the connected components of the space of mappings $S^1\to M$ 
that contain $\alpha$ and $\beta$, respectively.
Let $r: (S^1, *)\ \to (\no, \alpha)$ be a mapping with the adjoint $\hat r:S^1\times S^1\to M$, $\hat r_{|_{S^1\times 1}}=\alpha$. 

If 
$\ker \bigl(\hat r_*:\pi_1(S^1\times S^1)\to \pi_1(M)\bigr)\neq 1$, then using 
obstruction theory we get that $\hat r$ is homotopic to a mapping that passes through 
a mapping $S^1\to M$. This mapping of $S^1$ can be made disjoint from $\beta$ 
and $\mu_{1,0}(r,\beta)=0$ for 
such $r$. 

If $\hat r_*$ is injective then $r$ is homotopic 
to the  loop $\gamma_1^i \gamma_3^j\in (\no, \alpha)$, for some $i,j\in\Z$.  Here $\gamma_1$ is a self homotopy of $\alpha$ 
induced by one full rotation of the parameterizing circle and $\gamma_3$ is a 
self homotopy of $\alpha$ under which every point of $\alpha$ slides one full 
turn along the $S^1$ fiber of $F_g\times S^1\to F_g$ that contains the point, 
see the proof of~\cite[Lemma 6.11]{Chernovframed}. Clearly $\mu_{1,0}(\gamma_1, 
\beta)=0$ and $\mu_{1,0}(\gamma_3, \beta)$ equals to a $\pm 1$ times the class 
of the element of $\B$ that is obtained when $\alpha$ intersects $\beta$ under 
the deformation $\gamma_3$. Thus $\Im \mu_{1,0}=\Z\subset \bor_0(\B)$.

Similarly one shows that $\Im \mu_{0,1}=\Z\subset \bor_0(\B)$ and that in fact 
$\Im \mu_{1,0}=\Im \mu_{0,1}$. Thus $\Indet=\Z$ for this example. One can also 
show that $\bor_0(\B)=\bigoplus_{1}^{\infty}\Z$ in this case. 

Similar construction for $F_g=S^1\times S^1$ will give $\Indet=
\Z=\bor_0(\B)$.
\end{exs}

\begin{theodef}\label{mainthm} 
Let $\A$ be as above. 
Then there exists a function $\al\co  
\no\times \nd 
\setminus \Sigma \to \A$ such that
\begin{enumerate} 
\item[{\rm (a)}] $\al$ is constant on path connected components of 
$\no \times \nd \setminus \Sigma ;$ 
\item[{\rm (b)}] if $\gamma\co [a,b]\to \no\times \nd$ is a generic 
path 
such that 
$\gamma(a), \gamma(b)\not \in \Sigma$ and $t_i, i\in I,$ 
are 
the moments when 
$\gamma (t_i)\in \Sigma$ $($and hence $\gamma (t_i)\in 
\Sigma _0$ 
by the definition of the generic path$)$, then 
$$ 
\al(\gamma(b))- 
\al(\gamma(a))=q\Bigl (\sum _{i\in I}[\gamma 
(t_i)]\Bigr ) \in \A. 
$$ 
\end{enumerate} 

Moreover, these properties determine the function $\al$ 
uniquely up to an additive constant. We call such a function an {\rm affine 
linking invariant}. 
\end{theodef} 

We prove \theoref{mainthm} in Section~\ref{proofmainthm}. 

\begin{rem} After the first 
version~\cite{ChernovRudyakfirstversion} of this text appeared on the
electronic archive, U~Koschorke submitted the paper~\cite{Koschorke2}
where he used some invariants 
coming from the Hatcher--Quinn construction to study when two submanifolds of 
$M=S\times \R$ can be link homotoped to the disjoint $S$--levels.
\end{rem}

\begin{rem} It is rather easy to prove (see Corollary~\ref{imagealk})
that the mapping $\al\co \pi_0(\no \times \nd \setminus \Sigma) 
\to \A$ 
is always surjective. 
\end{rem}

\subsection{Relations to front propagation and winding numbers}\label 
{fronts} 
(a)\qua Let $STM$ denote the total space of the sphere tangent 
bundle over 
$M$, $\dim STM=2m-1$. In~\cite{ChernovRudyak} we defined 
the affine linking 
invariant 
$\pal$ for the mappings of $S^{m-1} \to STM$ 
that are homotopic to the inclusion of the fiber $S^{m-1}$ 
to $STM$. 
Because of the orientability of the bundle $STM \to M$, the 
homotopy class of 
this inclusion 
is invariant under the $\pi_1(M)$--action on $[S^{m-1}, STM] 
$.  
Using Theorem~\ref{pi0} we get that in this case $\bor_0(\B) 
=\Z$, and 
$\pal$ is 
exactly $\al$ 
for $N_1=N_2=S^{m-1}$ and the space $\no=\nd$ 
consisting of mappings $S^{m-1} \to STM$ as above. 

In~\cite{ChernovRudyak} we have shown that in this case 
$\eps(\Indet)=\Indet =0$ when $m$ is even or when $m$ is 
odd and $M$ is not a rational homology sphere. This shows that 
$\eps (\al)$ can indeed be $\Z$--valued in many cases where the mappings are 
not zero-homologous. 

This example is interesting since it is related to wave front propagation. Deep 
relations between front propagation and link theory were first discovered by 
V~Arnold~\cite{Arnold}, who observed that the isotopy type of the knot 
canonically associated to a front does not change under front propagation. 
A wave front on $M$ is a singular spherical hypersurface equipped with
a velocity vector 
field of the directions of the front propagation. The wave front represents the 
set of all points on $M$ that a certain information has reached at the given 
moment. Wave fronts on $M$ can be canonically lifted to $STM$ by mapping a point 
of the front to the point of $STM$ corresponding to the front velocity at this 
point. We assume that wave front propagation on $M$ is given by a time-dependent 
flow on $STM.$ (For example light propagation is indeed given by such a flow.)

Consider two wave fronts that originated at two different points on
$M$. The value of the 
$\al$ invariant on the canonical lifts to $STM$  of the two wave fronts equals 
to the algebraic number of times the front of the earlier event has passed 
through the birth point of the later event before the later event occurred, 
see~\cite{ChernovRudyak}. Thus if the $\al$ invariant is nonzero we conclude 
that the later event occurred after the earlier born front passed through its 
birth point. Philosophically this means that the later event has obtained the 
information about the earlier event carried by its wave front. Such events are 
called causally related. Thus the $\al$ invariant often allows one to detect 
that the events that created the wave fronts are causally related from the 
current picture of the wave front, without the knowledge of the time-dependent 
propagation law, times and places of events that created the wave fronts. 

First examples relating causality to the linking number were constructed by 
R~Low in the case of $M=\R^n$, see~\cite{Low1},~\cite{Low2},~\cite{Low3}. In 
this case the canonical lifts of the fronts are homologous into the end of 
$ST\R^n$. For such submanifolds the classical linking number can be defined as 
it was done by S~Tabachnikov~\cite{Tabachnikov} for $M=\R^3$. (This 
construction of the classical linking number was later generalized by 
U~Kaiser~\cite{Kaiserbook} to the case of arbitrary submanifolds of the linking 
dimension that are homologous into a boundary or into an end of the ambient 
manifold.) The modified Low conjecture says that two events are causally 
unrelated if and only if the lifts of their fronts are unlinked in the 
appropriate sense. For $M$ being a $2$--disk with holes strong results proving 
the conjecture for many cases were obtained by J~Natario and 
P~Tod~\cite{NatarioTod}.

(b)\qua The classical winding number of a curve around a point $p$ 
in $\R^2$ is the linking number between the curve and the 
$0$--cycle $\{p, \infty\}$ in $S^2=\R^2\cup \{\infty\}$. 
In~\cite{ChernovRudyakWinding} we considered the 
generalizations $\win (F, p)$ and $\wt \win (F, p)$ of the 
winding numbers of the mapping $F\co N_1^{m-1}\to M^m$ around a point $p\co \pt=N_2 
\to M$ to the case where $F_*([N_1])\neq 0\in H_*(M)$ and $M$ does not 
have ends that could play the role of the infinity. (The invariants 
$\wt \win$ and $\win$ are the generalizations of 
the winding number to the case where the observable point 
$p$ moves and is fixed in $M$, respectively.) 
We showed that 
affine 
winding numbers can be effectively used to estimate 
from below the number of times a wave front has passed 
through a point between two moments of time. 

One can verify that the generalized winding 
number $\wt \win$ also is included into the theory explored 
in this work, if we consider affine linking number for the 
case $N_2=\pt$. 
(It is clear, see~Theorem~\ref{pi0}, that in 
this case $\bor _0(\B)=\Z$.) 

It is easy to construct the version $\ov \al$ 
of the invariant $\al$ constructed in this paper, where 
$\ov \al$ 
will 
be a function on $\pi_0(\no \times \{*\}\setminus \Sigma)$ 
for 
some fixed mapping $*\in \nd$. 
A straightforward verification shows that 
$\ov \al$ is well-defined provided it takes values in the 
quotient group 
of 
$\bor 
_0 
(\B)$ by 
$\Im \bigl (\mu_{10} \co  \bor_1(\no)\otimes \bor_0(\nd)\to 
\bor _0 
(\B)\bigr )$. 
The invariant $\win$ constructed in~\cite{ChernovRudyakWinding} 
is a particular case of such $\ov \al$ where $N_2=\pt$.

\subsection[The invariant al is the universal
  Vassiliev--Goussarov invariant of order <=1]{The invariant $\al$ is the universal
  Vassiliev--Goussarov invariant of order $ \le 1$}
Fix a natural number $n$. Let 
$f=(f_1, f_2)\in \Sigma\subset \no 
\times \nd$ be such 
that 
$\Im (f_1)\cap \Im (f_2)$ consists of $n$ distinct 
double points of 
transverse intersection. 
Each double point can be resolved in two essentially 
different ways: 
positive 
and negative, where the sign is 
as in Definition~\ref{signalk}. Thus $f$ with $n$ such double 
points admits $2^n$ possible resolutions of the double 
points. A sign of the resolution 
is put to be $+$ if the number of negatively resolved 
double points is even, and it is put to be $-$ otherwise. 

Let $\Gamma$ be an abelian group. 
Let $\alpha$ be a $\Gamma$--valued {\it homotopy link invariant\/} of links from 
$\no\times \nd\setminus \Sigma$, ie a function $\alpha\co \pi_0(\no\times 
\nd\setminus \Sigma)\to \Gamma$. (Thus 
$\alpha$ does not change under homotopies of a linked pair that do not involve 
double points between different components.)
The {\it $n$th derivative $\alpha^{(n)}$ of $\alpha$\/} is a function on 
singular links 
with exactly $n$ distinct transverse intersection points between different 
component (and possibly many self-intersection points of the components). The 
value of $\alpha^{(n)}$ on such a singular link $f=(f_1, f_2)$ is defined as a 
sum (with appropriate signs) of the values of $\alpha$ on the nonsingular 
mappings obtained by the $2^n$ resolutions of double points.
The invariant $\alpha$ is said to be of {\em 
order $\le n-1$\/} (or \emph{Vassiliev--Goussarov invariant of order $\le 
n-1$}, \cite{Vassiliev,Goussarov1,Goussarov2}) if $\alpha^{(n)}$ is identically zero.

The invariant $\al$ is an $\A$--valued Vassiliev--Goussarov invariant of order 
$ \le 1$. To see this consider a singular link $f=(f_1, f_2)$ with exactly two 
transverse double points between different components. We denote by $f(\pm, 
\pm)$ the four nonsingular links obtained by resolutions of the two double 
points. We denote by $f(\cdot, \pm)$ the singular link with one transverse 
double point between different components obtained from $f$ by resolving the 
second singular point and keeping the first singular point intact.
To prove that $\al$ is an order $\leq 1$ invariants it suffices to show that 
$$\al(f(+,+))-\al(f(+, -))-\al(f(-,+))+\al(f(-,-))=0.$$
Rewrite this
expression as
\begin{multline*}
\bigl (\al(f(+,+))-\al(f(-,+))\bigr)-\bigl( \al(f(+,-))-\al(f(-,-))\bigr)\\
=\al^{(1)}(f(\cdot, +))-\al^{(1)}(f(\cdot, -)).
\end{multline*}
By the definition of 
$\al^{(1)}$ the values $\al^{(1)}(f(\cdot, +)), \al^{(1)}(f(\cdot, -))$ are 
equal to
$+1[f(\cdot, +)],$ and  $+1[f(\cdot, -)]\in \bor_0(\B)$. (Here
$[f(\cdot, +)], [f(\cdot, -)]\in \pi_0(\B)$ are the classes of the singular 
links with one transverse double point.) 
Clearly $[f(\cdot, +)]=[f(\cdot, -)]\in \pi_0(\B)$. Thus
\begin{multline*}
\al(f(+,+))-\al(f(+,-))-\al(f(-,+))+\al(f(-,-))\\
=\al^{(1)}(f(\cdot, +))-\al^{(1)}(f(\cdot, -))=0.
\end{multline*}
Furthermore, if $\alpha\co \pi_0(\no\times \nd \setminus 
\Sigma)\to 
\Gamma$ is a $\Gamma$--valued 
Vassiliev--Goussarov invariant of order $\le 1$, then the increment 
$\alpha^{(1)}(\gamma(t_0))=\Delta_{\alpha}(\gamma(t_0))$ of $\alpha$ under the 
positive crossing of $\Sigma$ at $\gamma (t_0)\in \Sigma_0$ 
depends only on 
the element of $\pi_0(\B)$ that corresponds to $\gamma 
(t_0)$. To see this consider a singular link $f$ with two double points obtained 
from $\gamma(t_0)$ by preserving the existing double point and creating a second 
double point by a homotopy. Since $\alpha$ is an order $\leq 1$ invariant
we have
\begin{align*}
0 &= \alpha(f(+,+))-\alpha(f(+,-))-\alpha(f(-,+))+\alpha(f(-,-)) \\
  &= \alpha^{(1)}(f(\cdot,+))-\alpha^{(1)}(f(\cdot, -)).
\end{align*}
Thus the increments into $\alpha$ of the 
positive crossings of the discriminant at $\alpha^{(1)}(f(\cdot, +))$ and at 
$\alpha^{(1)}(f(\cdot, -))$ are equal. Clearly we can change $\gamma (t_0)\in 
\Sigma_0$ to 
any other element in $\Sigma_0$ that is in the same component of $\pi_0(\B)$ by 
elementary homotopies 
that pass through an extra double point and by homotopies in $\pi_0(\B)$ that do 
not create extra double points between the two components. These operations do 
not change the value of $\alpha^{(1)}$ on an element in $\B$ and we see that 
$\alpha^{(1)}(\gamma(t_0))$ indeed depends only on the element of $\pi_0(\B)$ 
realized by $\gamma(t_0)$.

In particular, there exists 
a natural homomorphism $B\co \bor _0 (\B)\to \Gamma$ that 
sends the bordism class of 
$(+1)\gamma(t_0)\in \bor _0(\B)$ to $\alpha^{(1)}(\gamma(t_0))=\Delta_ 
{\alpha}(\gamma(t_0))$. Moreover, this homomorphism $B$ 
passes through 
the 
quotient homomorphism $q\co  \bor_0(\B) \to \A$ as in \ref{a}, 
and 
therefore we get 
a homomorphism $A\co  \A \to \Gamma$ with $A\circ q=B$, cf 
\theoref{alklkrelation} 
below. One verifies that 
$$ 
A (\al (f)-\al (f'))=\alpha(f)-\alpha(f'), 
$$ 
for all $f, f'\in (\no, \times \nd \setminus \Sigma)$. 

Clearly if $\alpha$ and $\alpha'$ are $\Gamma$--valued 
Vassiliev--Goussarov invariants of order $ \le 1$ such that 
$\alpha-\alpha'$ is a constant mapping, then the 
corresponding homomorphisms 
$A$ and $A'$ are equal. Let $f_0\in \no\times \nd \setminus 
\Sigma$ 
be a chosen preferred point. 
Then for every $f\in \no\times \nd \setminus \Sigma$ we 
have 
$\alpha(f)=\alpha(f_0)+A(\al (f))-A(\al (f_0))$. 
Thus $\al$ completely 
determines the values of $\alpha$ 
on all $f$ 
(modulo $\alpha(f_0)$), and hence $\al$ is a universal 
Vassiliev--Goussarov invariant of order $ \le 1$. 

In particular, $\al$ distinguishes all the elements 
$f, f'\in \no\times \nd\setminus \Sigma$ that can be 
distinguished via Vassiliev--Goussarov 
invariants of order $\le 1$ with values in an arbitrary 
group $\Gamma$.

\section{Proof of Theorem~\ref{mainthm}}\label{proofmainthm} 

\begin{defn}\label{sigma1} 
We define $\Sigma_1$ to be the subset (stratum) of $\Sigma$ 
consisting 
of all the pairs $(f_1, f_2)$ such that there exists 
precisely two 
pairs of points $y_1\in N_1, y_2\in N_2$ as in Definition~\ref{sigma0}. 
Here we assume 
that the two double points of the image are distinct. We also choose a 
base point $*$ of $\no\times \nd$ with $*\notin \Sigma$. 

Notice that $\Sigma_i, i=0,1,$ is the stratum of 
codimension $i$ 
in $\Sigma$. In particular, a generic path in $\no \times 
\nd$ intersects $\Sigma_0$ in a finite number of points, 
and a generic disk in $\no \times \nd$ intersects 
$\Sigma_1$ in a finite number of points. 
\end{defn} 

A generic path 
$\gamma \co [0,1]\rightarrow \no \times \nd$ 
that connects two points in $\no \times \nd \setminus 
\Sigma$ intersects $\Sigma_0$ in finitely many 
points $\gamma(t_j), 
j\in J,$ and all the intersection points are of the 
types described 
in~\ref{signalk}. Put 
\begin{equation} 
\Delta_{\al}(\gamma)=\sum_{j\in J}[\gamma(t_j)]\in \bor_0 
(\B). 
\end{equation} 
We let
\begin{align*}
A   &= \{(x,y)\in \R^2\bigm| x^2+y^2\le 1\}, &
B_1 &= \{(x,y)\in A \bigm| xy=0\},\\
B_2 &= \{(x,y)\in A \bigm| x=0\}, &
B_3 &= \{(0,0)\},\\
B_4 &= \emptyset. 
\end{align*}
We define a {\it regular disk} in $\no \times \nd$ 
to be a generically embedded disk $D$ such that 
$D\cap \Sigma=D\cap (\Sigma_0\cup\Sigma_1)$ 
and the triple $(D, D \cap \Sigma_0, D \cap \Sigma_1)$ 
is homeomorphic to a triple $(A, B, C), A\supset B \supset C$, where $B$
is one of $B_i$'s and $C$ is equal to $B_3$ or $B_4$.

\begin{lem}\label{small} 
Let $\beta$ be a generic loop that bounds a regular disk in 
$\no \times \nd $. Then $\Delta_{\al}(\beta)=0$. 
\end{lem} 

\begin{proof} It is easy to see that all the crossings of $\Sigma_0$ 
that happen along $\beta$ can be subdivided into pairs, 
such that 
the elements of $\pi_0 (\B)$ corresponding to the two 
crossings in a pair are equal and the signs of the 
corresponding crossings 
of $\Sigma _0$ are opposite. Hence the inputs into $\Delta 
_ 
{\al}(\beta)$ 
of the elements of $\bor _0(\B )$ corresponding 
to the two crossings in a pair cancel out 
and $\Delta_{\al}(\beta)=0$. 
\end{proof}

\begin{lem}\label{disk} 
Let $\beta$ be a generic loop that bounds a disk in 
$\no \times \nd$. Then $\Delta_{\al}(\beta)=0$. 
\end{lem} 

\begin{proof} Notice that the set $\Sigma\setminus(\Sigma_0 \cup \Sigma_1)$ is a 
subset of 
codimension $\ge 3$ in $\no\times \nd$. So, without loss of generality we can 
(using a small deformation of the disk) assume that the disk is the union 
of regular ones, 
cf Arnold~\cite{Arnoldcurves,Arnold}. Now the 
proof follows from Lemma \ref{small}. 
\end{proof} 

\begin{cor}\label{homomorphism} 
The invariant $\Delta_{\al}$ induces a well-defined 
homomorphism 
$$ 
\Delta_{\al}\co \pi_1 (\no \times \nd)=\pi_1 (\no \times \nd, 
*)\rightarrow \bor_0 
(\B). 
$$ 
\end{cor} 

\begin{proof} Since every element of $\pi_1(\no \times \nd, *)$ can 
be represented by a generic loop, the proof follows 
from Lemma \ref{disk}. 
\end{proof}

Now we give another description of $\Delta_{\al}$. Let 
$ 
\gl\co  \bor_1(\no \times \nd) \to \bor_0(B) 
$ 
be the homomorphism from \defref{a}. 

\begin{prop}\label{equality} 
The homomorphism $\Delta_{\al}\co \pi_1(\no \times \nd) \to 
\bor_0(\B)$ 
coincides with the homomorphism 
$$\pi_1(\no \times \nd) \stackrel{h}{\longrightarrow}
  \bor_1(\no \times \nd) \stackrel{\gl}{\longrightarrow}
  \bor_0(\B),$$ 
where $h$ is the Hurewicz homomorphism in the bordism 
theory $\bor_*(-)$ 
and $\gl$ is the homomorphism from Definition $\ref{a}$. 

\end{prop} 

\begin{proof} Given a generic loop $\ga$ in 
$(\no,*)$ and the 
constant loop $e$ in $(\nd,*)$, let $(\ga,e)$ be the 
corresponding loop in 
$(\no \times \nd,*)$. The homotopy class of $(\ga,e)$ gives 
us an element 
$[(\ga,e)]\in \pi_1 (\no \times \nd, *)$. Similarly, a 
generic loop $\gb$ 
in $(\nd,*)$ gives us an element $[(e,\gb)]\in \pi_1 (\no 
\times \nd, *)$. 

Because of the isomorphism $\pi_1(\no \times \nd)=\pi_1 
(\no) 
\times\pi_1(\nd)$, the classes $[(\alpha,e)]$ 
and $[(e,\beta)]$ generate the 
group 
$\pi_1(\no \times \nd, *)$. So, since $\Delta_{\al}$ is a 
homomorphism, 
it suffices to prove that 
$\Delta_{\al}[(\ga,e)]=(\gl\circ h)[(\ga,e)]$, and 
similarly for 
$[(e,\beta)]$. We do it here for the loops $(\ga, e)$ only.

We calculate $\Delta_{\al}[(\alpha, e)]\in \bor_0(\B)$. Fix 
a mapping 
$\ov e\co  N_2 \to M$ in $\nd$ and 
consider the mapping 
$$\overline \alpha\co S^1\times (N_1\sqcup N_2) \to M$$ 
such that $\overline \alpha \big|_{S^1\times N_1}=\alpha$ 
and $\overline \alpha \big |_{S^1\times N_2}$ 
coincides with the composition 
$$ 
\CD 
S^1\times N_2 @>\text{projection}>> N_2 @>\ov e>> M. 
\endCD 
$$ 
Without loss of generality we may assume that 
$\overline \alpha \big |_{S^1\times N_1}$ is transverse to 
$\ov e$. The inclusion $\no \to \no\times \nd,\, x\mapsto(x,e)$ allows 
us to regard $h(\ga)\in\bor_1(\no)$ as an element of $\bor_1(\no\times 
\nd)$. 

Now it is easy to see that 
\begin{equation}\label{twosetsequal} 
\Delta_{\al}[(\alpha, e)]= 
\mu_{10}\bigl (\bigl [ \overline \alpha \big|_{S^1\times 
N_1} 
\bigr ]\otimes \bigl [\ov e 
\bigr ]\bigr )=\gl(h(\ga))\in \bor_1(\B), 
\end{equation} 
where $\bigl[\overline \alpha \big |_{S^1\times N_1}\bigr]$ and $[\ov 
e]$ are the bordism classes of corresponding maps. 
\end{proof}

\begin{cor}\label{indet} 
$$ 
\Im\{\Delta_{\al}\co \pi_1 (\no \times \nd, *)\rightarrow 
\bor_0 
(\B)\}=\Indet \subset \bor_0(\B).\eqno{\qed} 
$$ 
\end{cor} 

Take an arbitrary point 
$f= (f^1_1,f^1_2)\in \no\times \nd \setminus \Sigma$ 
and choose a generic path $\gamma$ 
going from $*$ to $f$. We set 
\begin{equation*} 
\al(f)=q\bigl ( \Delta_{\al}(\gamma)\bigr )\in \A 
\end{equation*} 
where $q$ is the epimorphism from Definition~\ref{a}.

\begin{thm}\label{main} 
The function 
$$ 
\al\co \no\times \nd \setminus \Sigma\rightarrow \A 
$$ 
is an affine linking invariant. Furthermore, any other affine linking invariant 
$\wt{\al}$ coincides with $\al$ if $\wt{\al}(*)=0$. 
\end{thm} 

\begin{proof}
To show that $\al $ is constant on path components we must verify that the 
definition of $\al$ is independent on the choice of the path $\gamma$ that goes 
from $*$ to $f$. This is the same as to show that $q(\Delta_ 
{\al} (\varphi))=0$ for every closed generic loop $\varphi$ at $*$. But this 
follows from \corref{indet} directly. 

Furthermore, it is clear that $\al$ increases by $q([\gamma(t)]) \in\A $ under a 
transverse passage by a path $\gamma$ through the stratum $\Sigma_0$ at the 
point $\gamma(t)$. This yields property (b) of $\al$ from \theoref{mainthm}. The 
last claim is obvious. 
\end{proof}

\begin{rem}\label{affinerem} 
Clearly $\al$ depends on the choice of $*$. On the other hand, if we 
change $*,$ then new $\al$ and the old $\al$ will differ by an additive 
constant. This is the reason why we use the adjective ``affine''. 
\end{rem} 

Clearly, Theorem~\ref{mainthm} is a direct consequence 
of Theorem~\ref{main}.

\section[Relations between alk and the classical linking 
invariant]{Relations between $\al$ and the classical linking 
invariant $\lk$}\label{secrelation} 

Given a closed oriented manifold $N^n$ with the fundamental 
class $[N]\in 
H_n(M)$, we say that a map $f\co  N \to M$ is zero-homologous 
if 
$f_*([N])=0\in H_n(M)$. 

Let $\eps\co  \bor_0(\B) \to \Z$ be the homomorphism from 
\eqref{eps}. 

\begin{thm}\label{alklkrelation} 
Suppose that $\no$ and $\nd$ 
consist of zero-homologous mappings. Then $\eps(\Indet)=0$. 
Furthermore, 
for all $f=(f_1,f_2), f'=(f'_1,f'_2) \in \mathcal N=\no 
\times \nd\setminus 
\Sigma$, we have 
$$ 
\eps(\al (f)) - \eps(\al (f'))=\lk (f_1, f_2)- 
\lk (f'_1, 
f'_2)\in \Z. 
$$ 
\end{thm} 

\begin{proof}
Since $\no$ and $\nd$ consist of zero-homologous maps, the 
classical linking invariant $\lk\co  \no\times \nd \setminus 
\Sigma \to \Z$ 
is well-defined. Now, similarly to $\Delta_{\al}$, we 
define 
$$ 
\Delta_{\lk}\co  \pi_1(\no \times \nd,*) \to \Z, \quad 
\Delta_{\lk}(\gamma)=\sum_{i=1}^k\sigma(\gamma, t_i), 
$$ 
where the generic loop $\gamma$ in $(\no\times \nd, *)$ 
intersects 
$\Sigma_0\subset \Sigma \subset \no\times \nd$ in certain 
points $\gamma(t_1), 
\ldots, \gamma(t_k)$. (Here we use the notation $\gamma$ 
for the loop as well as 
for its homotopy class.) Since $\lk$ is well-defined, we 
conclude that 
$\Delta_{\lk}(\gamma)=0$ for all $\gamma$, ie $\sum\sigma 
(\gamma, 
t_i)=0$. 

Now, we have 
$$ 
\Delta_{\al}([\gamma])=\sum[\gamma(t_i)]\in \bor_0(\B). 
$$ 
So, in view of \eqref{eps} 
$$ 
\eps(\Delta_{\al}([\gamma]))=\sum\eps([\gamma(t_i)]) 
=\sum\sigma(\gamma, 
t_i)=0. 
$$ 
Thus, by Corollary~\ref{indet}, $$ 
\eps\left(\Indet\right)=\eps(\Im(\Delta_{\al}))=0. 
$$ 
Now take a generic path $\gamma$ which connects $f$ and $f'$. Then 
\begin{align*}
\eps (\al (f)){-}\eps(\al(f'))
  &=\eps\left(\sum[\gamma(t_i)]\right)=\sum\sigma(\gamma,t_i) \\
  &=\lk (f_1,f_2){-}\lk (f'_1,f'_2),
\end{align*}
which proves the second claim of the theorem.
\end{proof}

\begin{rems} Theorem~\ref{alklkrelation} demonstrates that, 
up to an 
additive constant, $\eps\circ\al$ is equal to the classical 
linking 
number $\lk$ whenever $\no$ and $\nd$ consist of zero-homologous 
mappings. 
So, $\al$ is an extension of the classical $\lk$--invariant 
of 
zero-homologous submanifolds. 

Since the homomorphism 
$\eps$ is the summation over the components of $\B$, we 
conclude that, for 
zero-homologous mappings, 
{\it $\al$ can be regarded as a splitting of the 
classical linking invariant} into 
a collection of independent invariants. 
In many cases it can be shown that $\bor _0(\B)$ is an 
infinitely generated 
abelian group (see Theorem~\ref{pi0}--Corollary~\ref{infinitePreissman}) 
and that the indeterminacy subgroup $\Indet$ is zero 
(see Section~\ref{Indet=0}). 
Since it is easy to show (see Corollary~\ref{imagealk}) that 
the mapping 
$\al\co \pi_0(\no\times \nd \setminus\Sigma)\rightarrow \A$ is 
always surjective, we see that for these cases the 
classical 
$\lk$ invariant of zero-homologous submanifolds splits into 
infinitely many independent invariants.

On the other hand, as it was explained in this paper, the invariant $\al$ 
exists regardless of whether the mappings are zero-homologous or not. 

Also, there are many cases where $\no, \nd$ do not
consist of zero-homologous mappings while nevertheless $\eps(\Indet) 
=0=\Indet$, and thus $\eps\circ\al$ is a $\Z$--valued invariant.
\end{rems}

\section{Examples where the indeterminacy subgroup vanishes} 
\label{Indet=0} 

Given the manifolds $M, N_1, N_2$ as in Section~\ref{affinesection}, we assume 
in 
addition that 
$N_1$ and $N_2$ are connected, and that $n_1n_2>0$.

\begin{thm}[Preissman]\label{hyper} 
Let $M$ be a closed manifold that admits a 
Riemannian metric of negative sectional curvature. 
Then the following holds: 

\begin{enumerate}
\item[\rm (i)] Every nontrivial abelian subgroup of $\pi_1(M)$ is an
infinite cyclic group.

\item[\rm (ii)] For every nontrivial abelian subgroup $A$ of 
$\pi$ there exists a unique 
abelian subgroup $B_A$ of $\pi$ which contains $A$ and is 
maximal with respect to 
this property. In fact, $B_A$ is the centralizer $Z(A)$ of 
$A$ in $\pi_1(M)$. 
\end{enumerate}

\end{thm} 

\begin{proof} See do Carmo~\cite{docarmo} or the original paper by 
Preissman~\cite{Preissman}. 
\end{proof} 

\begin{defn}\label{Preissmangroup} A finitely generated group $\pi$ is called 
a {\it Preissman group} if it satisfies the properties (i) 
and 
(ii) from \theoref{hyper}. 
\end{defn} 

\begin{prop}\label{commut} 
Let $\pi$ be a Preissman group. Let $\ga, \gb\in 
\pi$ be such that 
$\ga\gb\ne \gb\ga$, and let $\gamma\in \pi$ be such that 
$\ga\gamma=\gamma\ga$ and $\gb\gamma=\gamma\gb$. Then 
$\gamma=e$. 
\end{prop} 

\begin{proof} Suppose that $\gamma \ne e$. Let $G=\{x\}$ be the 
(unique) maximal 
cyclic subgroup of $\pi$ which contains $\gamma$. Since 
$\ga\gamma=\gamma\ga$, the 
subgroup $\{\ga, \gamma\}$ of $\pi$ is contained in $G$, 
and so 
$\ga=x^m$ 
for some $m$. Similarly, $\gb=x^k$, and thus 
$\ga\gb=\gb\ga$. This 
is a 
contradiction. 
\end{proof}

\begin{thm}\label{ZeroIndet} 
For $M, N_1, N_2$ as above, suppose that $\pi_1(M)$ is 
a Preissman group, and that $\pi_i(M)=0$ 
for $2\le i \le 
1+\max\{n_1, n_2\}$. Then the indeterminacy subgroup 
$\Indet\subset 
\bor_0(\B)$ is the zero subgroup, $\Indet =\{0\}\subset 
\bor_0(\B)$. 
\end{thm} 

\begin{proof} Throughout the proof we denote $\pi_1(M)$ by $\pi$. We 
must prove that, for 
every $\ga\in 
\pi_1(\no)$ and $\gb\in \pi_1 (\nd)$, 
\begin{equation}\label{delta=0} 
\Delta_{\al}[(\alpha, e)]=0=\Delta_{\al}[(e,\beta)]
\end{equation} 
(see \lemref{indet}). 

We prove the first equality from \eqref{delta=0} only, 
the second equality can be proved in the similar way. Fix 
$\phi_i\in\NN_i$ for $i=1,2$, and consider a loop $\ga$ in
$(\no,\phi_1)$. Let $\wt \ga\co S^1\times N_1 \to M$ be the adjoint map,
$\wt \ga(t,n)=\ga(t (n)$. Since $\pi_i(M)=0$ for $2\le i 
\le 1+n_1$, it follows from the elementary obstruction 
theory that the 
homomorphism 
$$\wt \ga_*\co  \pi_1(S^1\times N) \to \pi$$ 
completely determines the homotopy class of $\wt \ga$. We 
use the 
isomorphism $\pi_1(S^1\times N) \simeq \pi_1(S^1)\times 
\pi_1(N)$ and set 
\begin{equation}\label{gamma} 
\gamma=\wt \ga_*(\iota, e), 
\end{equation} 
where $\iota\in \pi_1(S^1)$ is the generator. 

\begin{lem}\label{gamma=0} 
If $\gamma=e$, then $\Delta_{\al}[(\ga,e)]=0$. 
\end{lem} 

\begin{proof} Indeed, in this case $\wt \ga\co  S^1\times N_1 \to M$ is 
homotopic to the map 
$$ 
\CD 
S^1 \times N_1 @>\text{proj}>> N_1 @>\phi_1 >> M, 
\endCD 
$$ 
because both maps induce the same homomorphism of 
fundamental groups. So, since $n_1+n_2=m-1<m$, 
there is a generic map $\wh\ga$ homotopic to $\wt \ga$ such 
that $\wh\ga(S^1\times 
N_1)$ does not meet $\phi_2(N_2)$. 
\end{proof} 

\begin{lem}\label{im=z} 
If $\Im \wt \ga_*=\Z \subset \pi$, then $\Delta_{\al} 
[(\ga,e)]=0$. 
\end{lem} 

\begin{proof} If $\Im \wt \ga_*=\Z$, then $\wt \ga_*$ can be 
decomposed as 
\begin{equation}\label{z} 
\pi_1(S^1\times N_1) \to \Z \subset \pi. 
\end{equation} 
Since $S^1=K(\Z,1)$, the map $\pi_1(S^1\times N_1) \to \Z$ 
in \eqref 
{z} can be induced by a map $\gf\co S^1\times N_1 \to S^1$. 
Furthermore, 
the 
inclusion $\Z\subset \pi$ in \eqref{z} can be induced by a 
map 
(inclusion) 
$\psi\co  
S^1 \to M$, and we can assume that $\psi(S^1)$ does not 
meet 
$\phi_2(N_2)$ since $m-n_2>1$. Now, the map 
$\wt \ga\co  S^1 \times N_1 
\to M$ is homotopic to a map 
\[ 
\wh \ga\co  S^1\times N_1 \xrightarrow{\gf} S^1 \xrightarrow 
{\psi} M. 
\] 
Clearly, $\wh \ga$ does not meet 
$\phi_2(N_2)$ and, moreover, any small perturbation of $\wh 
\ga$ does 
not. 
Thus, $\Delta_{\al}[(\ga,e)]=0$. 
\end{proof} 

Now we finish the proof of \theoref{ZeroIndet}. 
Consider the homomorphism 
$$ 
(\phi_1)_*\co  \pi_1(N_1) \to \pi. 
$$
If $\Im (\phi_1)_*$ is non-abelian, then $\gamma=e\in\pi$ 
by \propref{commut}, because 
$\gamma$ commutes with $\Im (\phi_1)_*$. So, $\Delta_{\al} 
[(\ga,e)]=0$ by 
\lemref{gamma=0}. 

Furthermore, assume that $\Im (\phi_1)_*$ is abelian. Since 
$\gamma$ commutes 
with $\Im (\phi_1)_*$ we conclude that $\Im(\wt \ga_*)$ is 
an abelian subgroup 
of $\pi$. So, $\Im(\wt \ga_*)=\Z$ or $0$ because $\pi$ is 
a Preissman group. If $\Im(\wt \ga_*)=\Z$ then 
$\Delta_{\al}[(\ga,e)]=0$ (by \lemref{im=z}), and if $\Im(\wt \ga_*)=0$, 
then $\gamma=e\in\pi$ and $\Delta_{\al}[(\ga,e)]=0$ 
(by~\lemref{gamma=0}). 
\end{proof} 

\begin{rem} Assume that $n_1+n_2+1\neq m$ but otherwise the manifolds $N_1, N_2, 
M$ satisfy all the conditions of Theorem~\ref{ZeroIndet}.
Assume moreover that $\max(n_1, n_2)< m-1$ and that there exists $(f_1, f_2)\in 
\no\times \nd$ with $\Im f_1\cap \Im f_2=\emptyset$. Then the straightforward 
modifications of the proof of the Theorems~\ref{ZeroIndet} shows that $\Im 
\mu_{1,0}=\Im \mu_{0,1}=0\subset \bor_{n_1+n_2+1+0-m}(\B)$.
\end{rem}

\begin{thm}\label{zeroindet2} 
Let $M$, $N_1$ and $N_2$ be as above with an extra condition that $\pi_i(M)$ is 
finite for $i=1, \dots, \max (n_1, n_2)+1$, then the indeterminacy subgroup 
$\Indet\subset \bor_0(\B)$ is the zero subgroup. 
\end{thm} 
\begin{proof} Elementary obstruction theory implies that a mapping 
$N_i\times*\to M$ (for $i=1,2$) has are finitely many non-homotopic
extensions to a
mapping $N_i\times S^1\to M$ (for $i=1,2$). Thus 
$\pi_1(\no\times \nd)$ is finite and hence
$$\Im\bigl(\Delta_{\al}\co \pi_1(\no\times \nd, *)\to
  \bor_0(\B)\bigr)=\Indet\subset \bor_0(\B)$$
is a finite subgroup. Since $\bor_0(\B)$ is torsion free, we get 
that $\Indet=0$. 
\end{proof} 

\begin{rem} Assume that $n_1+n_2+1\neq m$ but otherwise the manifolds $N_1, N_2, 
M$ satisfy all the conditions of Theorem~\ref{zeroindet2}. Then the 
straightforward modification of the proof of Theorem~\ref{zeroindet2} shows that 
$\Im \mu_{1,0}=\Im \mu_{0,1}=0\subset \bor_{n_1+n_2+1+0-m}(\B)$.
\end{rem}

\begin{exs} 
(1)\qua It is well-known that 
every closed manifold $M^m$ that admits a complete 
Riemannian metric of negative sectional curvature 
has the universal covering space homeomorphic to $\R^m$, and thus 
$\pi_i(M)=0$, for $i>1$, while $\pi_1(M)$ is a Preissman group. 

Combining this with 
Theorems~\ref{hyper} and~\ref{ZeroIndet} we get that for 
such $M$ the group 
$\Indet=0$ for all $N_1, N_2, \no, \nd$, 
and hence $\A$ is a free abelian group. 

To construct more manifolds with Preissman fundamental 
group, notice 
that the total space of a locally trivial bundle $F\to E 
\to B$ has the 
Preismann fundamental group if $F$ is simply-connected and 
$\pi_1(B)$ is 
Preissman.

(2)\qua Let $M^m$ be an oriented base space of a finite covering
map $\Sigma^m\to M$, where $\Sigma^m$ is a simply connected
$m$--dimensional rational homology sphere.  The Serre--Hurewicz
theorem implies that $\pi_k(\Sigma^m)$ is finite, for $k=2, \ldots,
m-1$, and therefore $\pi_k(M)$ is finite, for $k=1, \ldots, m-1$.
Now Theorem~\ref{zeroindet2} implies that $\Indet=0$ for all $N_1, N_2,
\no, \nd$, provided that $n_1, n_2>0$.

Theorem~\ref{alklkrelation} implies that in all these 
cases, if $\no, \nd$ consist of zero-homologous mappings, 
then 
$\al$ is a 
splitting of the classical linking invariant $\lk$ into a 
direct sum of independent $\Z$--valued invariants. 
\corref{infinitePreissman} says that this 
direct sum is 
infinite when 
$M$ is closed and admits a complete 
metric of negative sectional curvature; the images of 
$\pi_1(N_1), \pi_1(N_2)$ in $\pi_1(M)$ 
under the homomorphisms induced by the mappings $f_1\in \no, 
f_2\in \nd$ are nontrivial; and $\big |H_1(M)/\bigl (\Im f_{1*}(H_1(N_1))+\Im 
f_{2*}(H_1(N_2))\bigr)\big |=\infty.$ 

\end{exs}

\section[Description of pi_0(B) and omega_0(B)]{Description of $\pi_0(\B)$ and $\bor_0(\B)$}\label 
{B}\label{secdescription} 

In Section~\ref{secdescription} we assume $N_1$ and $N_2$ are connected;
\emph{however we do not assume that $n_1+n_2+1=m$}.  The goal of this
section is to give a method of evaluation of $\pi_0(\B)$ (and therefore of
$\bor_0(\B)$ since the last group is the free abelian group with the base
$\pi_0(\B)$) and, in particular, to demonstrate that $\pi_0(\B)$ can be
infinite in spite of $\pi_0(\no)$ and $\pi_0(\nd)$ being one-point sets.

We must recall some facts from elementary homotopy 
theory. 
Given two pointed spaces $X$ and $Y$, let $[X,Y]\bul$ and 
$[X,Y]$ 
be the set of pointed homotopy classes of pointed maps $X 
\to Y$ and 
the set of 
unpointed homotopy classes, respectively. Then 
$\pi_1(Y)$ acts 
on $[X,Y]\bul$ in a usual way, and 
\begin{equation}\label{quot} 
[X,Y]\bul/\pi_1(Y)=[X,Y] 
\end{equation} 
(see \cite{Spanier}, for example). 

Furthermore, 
$\pi_1(X)$ acts on $[X,X]$, and therefore we get a right 
$\pi_1(X)$--action on $[X,Y]$ via the composition map $[X,X] 
\times [X,Y]\to [X,Y]$. In greater detail, we apply 
$\ga\in \pi_1(X)$ 
to the 
homotopy class $1_X\in [X,X]$ and get a map (homotopy class) 
$\ga(1_X)$. Now, given $f\in [X,Y]$, we define $f\ga=(f)\ga\in 
[X,Y]$ to 
be the composition $f\circ\ga(1_X)$. It is worthy to 
mention that, 
for the standard 
$\pi_1(Y)$--action on $[X,Y]$ we have 
\begin{equation}\label{inverse} 
\gb(f)=\gb(1_Y)\circ f, \quad \gb \in \pi_1(Y), f\in 
[X,Y]. 
\end{equation} 
Furthermore, given $\ga \in \pi_1(X)$ and $f\in [X,Y]\bul$, 
we have 
\begin{equation}\label{orbit} 
f\ga=f_*(\ga)f 
\end{equation} 
and so any orbit of the (right) $\pi_1(X)$--action is contained 
in an orbit of the $\pi_1(Y)$--action. 

\begin{lemma}\label{actionscommute} 
The $\pi_1(Y)$--action on $[X,Y]\bul$ and the $\pi_1(X)$--action 
on $[X,Y]\bul$ commute. 
\end{lemma} 

\proof Choose $\ga\in \pi_1(X), \gb \in \pi_1(Y)$ and consider 
the 
composition 
$$ 
\CD 
X @>\ga(1_X)>> X @>f>> Y @>\gb(1_Y)>> Y. 
\endCD 
$$ 
Then, in view of \eqref{inverse}, we have
$$(\gb(f))\ga = 
( \gb(1_Y)\circ f)\circ\ga(1_X) 
=\gb(1_Y)\circ(f\circ\ga(1_X)) 
=\gb((f)\ga).\eqno{\qed}$$

Now, consider two pointed spaces $X_1$ and $X_2$ and the 
actions 
\begin{equation*} 
[X_1,Y]\bul\times\pi_1(X_1)\to [X_1,Y]\bul \quad \text{ 
and } \quad 
[X_2,Y]\bul\times\pi_1(X_2)\to [X_2,Y]\bul . 
\end{equation*} 
Together these actions yield the action 
\begin{equation}\label{direct} 
( [X_1,Y]\bul\times 
[X_2,Y]\bul)\times (\pi_1(X_1)\times \pi_1(X_2))\to 
[X_1,Y]\bul\times[X_2,Y]\bul . 
\end{equation} 
Because of \lemref{actionscommute}, this action commutes 
with the diagonal 
$\pi_1(Y)$--action on $[X_1,Y]\bul\times [X_2,Y]\bul$. 

\begin{rem} 
Since $[X_1,Y ]\bul \times [X_2,Y]\bul=[X_1\vee X_2,Y] 
\bul$, we have 
the action of the group $\pi_1(X_1\vee X_2)=\pi_1(X_1)*\pi_1 
(X_2)$ on 
$[X_1,Y]\bul \times [X_2, Y]\bul$. 
Let $i\co  \pi_1(X_1) \subset \pi_1(X_1)*\pi_1(X_2)$ be the 
standard 
inclusion. 
We emphasize that the $i(\pi_1(X_1))$--action 
on $[X_1,Y]\bul \times [X_2, Y]\bul$ does {\it not} 
coincide with 
the action coming from \eqref{direct}. 
\end{rem} 

Given $f\in [X,Y]\bul$, let 
\begin{equation}\label{stabilizer} 
S_f=\{\ga\in \pi_1(Y)\bigm|\ga f=f\}. 
\end{equation} 
\begin{prop}\label{stationary} 
$S_f\subset Z(f_*(\pi_1(X))$, where $ZG$ denotes the centralizer 
of $G$ in 
$\pi_1(Y)$. 
\end{prop} 

\begin{proof} Let $H$ be the set of homomorphisms $\pi_1(X) \to \pi_1 
(Y)$ induced 
by pointed maps $X \to Y$. The $\pi_1(Y)$--action on $[X,Y] 
\bul$ induces 
the $\pi_1(Y)$--action on $H$. Here for $\ga\in\pi_1(Y)$ we 
have 
$$ 
(\ga(f_*))(u)=\ga(f_*(u))\ga^{-1}, 
$$ 
where $u\in \pi_1(X)$ and $f\co  X \to Y$. Thus, if $\ga\in 
S_f$ then 
$(\ga(f_*))(u)=\ga(f_*(u))\ga^{-1}=f_*(u)$, ie $\ga \in 
Z(f_*(\pi_1(X))$. 
\end{proof}

Choose base points in $N_1, N_2$ and $M$ and let $\mathcal 
N_i\bul , i=1,2,$ be the space of all pointed maps $f\co  N_i 
\to M$ such 
that the 
unpointed map $f$ belongs to $\NN_i$. 
Clearly, the subset $\pi_0(\NN_i\bul)$ of $[N_i, M]\bul$ is 
invariant 
with 
respect to $\pi_1(N_i)$-- and $\pi_1(M)$--actions. Moreover, 
because of 
\eqref{quot}, $\pi_0(\NN_i\bul)$ is the orbit of $\pi_1(M)$--action on 
$[N_i, M]\bul$. 

Clearly, the bijection $[N_1\vee N_2, M]\bul=[N_1,M]\bul 
\times 
[N_2, M]\bul$ 
converts the $\pi_1(M)$--action on $[N_1\vee N_2,M]$ into 
the diagonal 
action 
on $[N_1,M]\bul \times [N_2, M]\bul$. 

\begin{thm}\label{pi0} 
$$ 
\pi_0(\B)= 
(\pi_1(N_1)\times\pi_1(N_2))\big\backslash\left(\pi_0(\no\bul)\times 
\pi_0(\nd\bul)\right)\big\slash \pi_1(M). 
$$ 
\end{thm} 

\begin{proof} Consider the inclusion $j\co \no\bul\times \nd\bul \to 
\B$. Since 
$M$ and $N_i$ are connected, the 
map $j_*\co  \pi_0(\no\bul\times \nd\bul) \to \pi_0(\B)$ is 
surjective. 
Now, consider two points $b,b'\in \no\bul \times \nd\bul$ 
and set 
$j(b)=(\rho_1,\rho_2, \phi_1, \phi_2)$ 
and $j(b')=(\rho'_1,\rho'_2, \phi'_1, \phi'_2)\in \B$. Then 
$\phi_1\rho_1=\phi'_1\rho'_1$. Suppose that there is 
a path from $b$ to $b'$ in $\B$. Consider this path as a
quadruple $(\rho_1(t),\rho_2(t), \phi_1(t), \phi_2(t)), 
t\in I$. 
Now, the loop $\rho_i(t)$ gives us the loop (homotopy 
class) 
$\ga\in \pi_1(N_i), i=1,2$, while the loop 
$$ 
(\phi_1(t))(\rho_1(t))=(\phi_2(t))(\rho_2(t)) 
$$ 
gives us an element $\gb\in \pi_1(M)$. Clearly, 
$\beta([b])(\ga_1,\ga_2)=[b']$ 
where $[b]$ and $[b']$ are components of $\no\bul \times 
\nd\bul$ containing 
$b$ and $b'$, respectively. 
\end{proof} 

\begin{cor}\label{imagealk} 
The mapping $\al\co \pi_0(\no \times \nd\setminus \Sigma)\to 
\A$ of Theorem~$\ref{mainthm}$ is surjective. 
\end{cor} 

\begin{proof} Take an arbitrary point $b\in \B$ and consider the class $[b]\in 
\bor_0(\B)$. Clearly it suffices to show that, given $(f_1, f_2)\in 
\no\times \nd \setminus \Sigma$, 
there exist generic paths 
$\gamma_i \co I\to \no\times \nd , i=1,2,$ starting at 
$(f_1, f_2)$ 
such that $\gamma_i(t)$ intersects $\Sigma_0$ for exactly one 
value $t_i\in (0,1), i=1,2,$ and 
$[\gamma_1(t_1)]=[b], [\gamma_2(t_2)]=-[b]$. Clearly there exists 
$b'=[\phi_1, \phi_2, \rho _1, \rho _2]\in \B$ such that 
$\pm [b]'=[b]$ and that 
one of the two possible resolutions of the double point 
between $\phi_1$ and $\phi_2$ is isotopic to $(f_1, f_2)$. 
Let $(\phi'_1, \phi'_2)$ be this resolution and 
let $\tilde \gamma$ be the isotopy path that connects 
$(f_1, f_2)$ to $(\phi'_1, \phi'_2)$. 

\begin{figure}[ht!]\small\anchor{change.fig}
\begin{center} 
\psfraga <-2pt,0pt> {f1}{$\phi'_1$}
\psfraga <0pt,2pt> {f2}{$\phi'_2$}
\includegraphics[width=10cm]{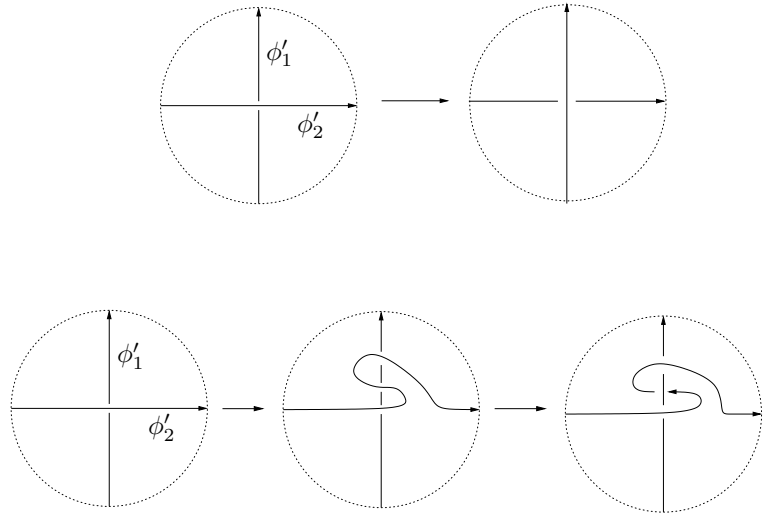} 
\end{center} 
\caption{The paths $\tilde \gamma_+$ and $\tilde \gamma_-$ 
in the $1$--dimensional case}\label{change.fig} 
\end{figure}

It is easy to construct (see \figref{change.fig}) a 
local 
deformation $\tilde \gamma_+$ (respectively, $\tilde \gamma_-$) of 
$(\phi'_1, \phi'_2)$ that involves only one 
passage through a double point and such that the sign of 
the corresponding crossing of $\Sigma_0$ is positive 
(respectively, negative) and the associated element of $\pi_0(\B)$ 
is the same as the one for $b'$. 
The desired paths $\gamma_i$ are the paths $\tilde 
\gamma_+\tilde \gamma$ and 
$\tilde \gamma_-\tilde \gamma$. 
\end{proof} 

\begin{cor}\label{infin} 
Assume that $N_1$ and $N_2$ are simply connected and that 
$\pi_1(M)$ 
is infinite. Also, assume that, for both $i=1$ and $i=2$, 
the 
$\pi_1(M)$--action 
on $\pi_0(\NN_i\bul)$ is free. Then $\pi_0(B)$ is an 
infinite set. \qed
\end{cor} 

\begin{ex} Let $N_1=N_2=S^n$, and let $M$ be a manifold 
which has 
the 
$(n+2)$--type of the wedge $S^1\vee S^n$. Fix a map $f\co S^n 
\to M$ which 
induces an isomorphism in $n$--dimensional homology. Let 
$\NN=\no=\nd$ 
be the path component of $f$ in the space of smooth maps 
$S^n \to M$. 
Then $\pi_0(B)$ is an infinite set because of \corref 
{infin}. 
\end{ex} 

Another series of cases with $\pi_0(\B)$ infinite 
appears as follows.

Take $f_i\in \pi_0(\NN_i\bul), i=1,2,$ and let $G_i$ be the 
subgroup of 
$\pi_1(M)$ generated by $S_{f_i}$ and $f_{i*}(\pi_1(N_i))$, 
where 
$S_{f_i}$ is the stabilizer of $f_i$ as in \eqref 
{stabilizer}.

\begin{thm}\label{doublecoset} 
If the Reidemeister set of double cosets 
$$ 
G_2\backslash\pi_1(M)\slash G_1 
$$ 
is infinite, then $\pi_0(\B)$ is an infinite set and $\bor_0 
(\B)$ 
is an infinitely 
generated abelian group. 
\end{thm} 

\begin{proof}Take any $u,v \in \pi_1(M)$ such that the elements
$(uf_1,f_2)$ and  
$(vf_1,f_2)$ of $\pi_0(\no\bul)\times \pi_0(\nd\bul)$ give us the same element 
of $\pi_0(\B)$.

It 
follows from \theoref{pi0} 
that there exists $(\ga,\gb, \gamma)\in\pi_1(N_1)\times 
\pi_1(N_2)\times 
\pi_1(M)$ such that 
$\gamma(uf_1,f_2)(\ga,\gb)=(vf_1,f_2)$ or 
\begin{equation*} 
\gamma uf_1\ga = vf_1 \text{ in }\pi_0(\no\bul) \quad \text 
{ and } 
\quad \gamma f_2\gb = f_2 \text{ in }\pi_0(\nd\bul). 
\end{equation*} 
Because of \lemref{actionscommute} and \eqref{orbit}, 
$\gamma 
f_2\gb=\gamma f_{2*}(\gb) f_2$. So, $\gamma f_ 
{2*}(\gb) \in S_{f_2}$, and therefore $\gamma\in G_2$. 

Similarly, $v^{-1}\gamma u f_1\ga=f_1$, ie $v^{-1} 
\gamma u 
f_{1*}(\ga)f_1=f_1$, and hence $v^{-1}\gamma u f_{1*}(\ga) 
\in S_{f_1}$. 
Thus $v^{-1}\gamma u \in G_1$, ie $\gamma u \in v G_1$, 
and $G_2 u 
\cap v G_1\ne \varnothing$. 

This means that the images of $u$ and $v$ in 
$G_2\backslash\pi_1(M)\slash G_1$ coincide. This completes 
the proof. 
\end{proof} 

\begin{cor} 
Let $h\co \pi_1(M) \to H_1(M)$ be the Hurewicz homomorphism. 
Then 
$\pi_0(B)$ is infinite whenever the group $H_1(M)/(h(G_1)+h 
(G_2))$ is 
infinite.
\end{cor} 

\begin{cor}\label{infinitePreissman} 
Assume that $\pi_1(M)$ is a Preissman 
group. Let $\no\times \nd$ be such that 
for $(f_1, f_2)\in \no \times \nd$ the subgroups $f_{1*} 
(\pi_1(N_1)), f_{2*}(\pi_1(N_2))\subset \pi_1(M)$ 
are both nontrivial 
and such that $H_1(M)\slash \bigl( f_{1*}(H_1(N_1))+f_{2*} 
(H_1(N_2))\bigr)$ is infinite. Then $\pi_0(\B)$ 
is infinite and hence $\bor_0(\B)$ is an infinitely 
generated abelian group. 
\end{cor} 

\begin{proof} It suffices to prove that each of the groups $f_{i*}(\pi_1 
(\NN_i)), 
i=1,2,$ is a subgroup of finite index in $G_i$. For sake of 
simplicity, we fix $i$ and denote $G_i$ and $f_{i*}(\pi_1 
(\NN_i))$ by $G$ 
and $H$, respectively, and prove that the index $[G:H]$ is 
finite. 

\textbf{Case 1}\qua Assume that the group $H$ is abelian. 
Then since 
$\pi_1(M)$ is a Preissman group, we conclude that the 
centralizer $ZH$ 
of $H$ is a cyclic subgroup which contains $H$. Since 
$H\ne\{e\}$, $H$ is 
a subgroup of finite index in $ZH$. Thus, we are done 
since, by 
\propref{stationary}, $S_f\subset ZH$ and so $G\subset 
ZH$. 

\textbf{Case 2}\qua Assume that the group $H$ is non-abelian. 
Then, by \propref{commut}, $ZH=\{e\}$. Thus, by \propref 
{stationary}, 
$S_f=\{e\}$, ie $G=H$. 
\end{proof} 

\section{The pairing $\mu_{i,j}$ as an obstruction for double point 
elimination and as a splitting of the intersection 
pairing}\label{disjointmappings} 

In this section we do not assume that $\dim N_1^{n_1}+\dim N_2^{n_2}=\dim M^m$. 

A classical question is: given two $C^{\infty}$ maps $f_1\co  N_1 \to M$
and $f_2\co  
N_2 
\to M$, can we make their images disjoint by homotopies of $f_1$ and 
$f_2$? 
The necessary condition for the existence of such homotopies is that the 
intersection of the homology classes 
$f_{1*}([N_1])$ and $f_{2*}([N_2])$ is zero in $H_{m-n_1- 
n_2}(M)$. It is well-known that this condition is 
not sufficient in general. 

From now on we assume that $n_1+n_2=m$ and consider the pairing 
$$ 
\mu=\mu_{0,0}\co  \bor_0(\no)\otimes\bor_0(\nd) \to \bor_{n_1+n_2+0+0-m} 
(\B)=\bor_0(\B). 
$$
The mapping $e\co  \B \to M$ that maps $(\phi_1, 
\phi_2, \rho_1, \rho_2)$ to 
$\phi_1\circ \rho_1(\pt)=\phi_2\circ \rho_2 (\pt)\in M$ 
induces the (augmentation) map 
$\eps \co \bor _*(\B )\to \bor _*(M)$. 
Clearly $\tau \circ \eps (\mu([f_1], [f_2]))$ is equal to 
the intersection of 
$f_{1*}([N_1])$ and $f_{2*}([N_2])$ in $H_{m-n_1-n_2}(M)$, 
where $\tau$ is the Steenrod--Thom homomorphism. 

We observe that the necessary condition 
for the existence of $f_1\in \no, f_2\in \nd$ with $\Im 
(f_1)\cap \Im (f_2)=\emptyset$ is that 
$$ 
\mu([\bar f_1], [\bar f_2])=0\in \bor _0(\B), 
$$ 
for some (and therefore for all) 
$\bar f_1\in \no, \bar f_2\in \nd$. Clearly this condition 
is much stronger than the vanishing of intersection 
condition. (Related results for intersections of mappings 
of $S^1$ into an oriented $2$--dimensional surface were previously obtained 
by Turaev and Viro~\cite{Turaevloops,TuraevViro}.) 

In fact, the following Theorem~\ref{whitneytrick} describes a big class of 
manifolds $M$ where the 
condition 
$$ 
\mu([\bar f_1], [\bar f_2])=0\in \bor _0(\B) 
$$ 
is necessary and sufficient for the existence of $f_1\in \no$ and 
$f_2\in \nd$ with disjoint images. 

When $N_1$ and $N_2$ are simply connected and $n_1+n_2=m$ such $f_1\in \no, 
f_2\in \nd$ with disjoint images do exist if and only if the 
$\Z[\pi_1(M)]$--valued intersection index of the lifts of the mappings $\bar f_1$ 
and $\bar f_2$ to the universal covering $\wt M\to M$ vanishes (see 
Kervaire~\cite{Kervaire}, for example). In the conditions of
Theorem~\ref{whitneytrick} $N_1$ and $N_2$ are not assumed to be simply
connected. Thus the mappings $\bar f_1$ and $\bar f_2$ can not be lifted
to the universal covering and Kervaire's approach to this classical
problem does not work.

Observe also that if $\no, \nd$ consist of mappings homotopic to immersions, 
then the statement 
of Theorem~\ref{whitneytrick} follows immediately from Theorem 2.2 of the 
work~\cite{HQ} of 
Hatcher and Quinn (see Section~\ref{hq}).

\begin{thm}\label{whitneytrick} 
Let $N_1^{n_1}, N_2^{n_2}$, $n_1, n_2>2$, be closed manifolds and let $M^m$ be 
a {\rm(}not necessarily closed\/{\rm)} 
manifold such that $n_1+n_2=m$, $\pi_1(M)$ is Preissman
{\rm(}see Definition~\ref{Preissmangroup}\/{\rm)} and $\pi_i(M)=0$ for $i=2, \dots, \max(n_1,
n_2)+1$.
Let $\no, \nd$ be connected components of the space of 
mappings of $N_1$ and of $N_2$ into $M$. Then $\mu\co \bor_0(\no)\otimes 
\bor_0(\nd)\to \bor_0(\B)=\Z$ is the zero pairing if 
and only there exists $f_1\in \no$, $f_2\in \nd$ with $\Im f_1\cap \Im 
f_2=\emptyset$. 
\end{thm} 

\begin{rem} In particular, by the Preissman and Hadamard Theorems, 
see Theorem~\ref{hyper} and~\cite{docarmo}, $\pi_1(M)$ is Preissman
and $\pi_j(M)=0$, $j>1$ for all closed $M$ that admit a Riemannian metric
of negative sectional curvature; and any such $M$ satisfies all the
conditions of Theorem~\ref{whitneytrick}. 
\end{rem} 

\begin{ex} Here we give an example where neither Kervaire non
Hatcher--Quinn approaches work while \theoref{whitneytrick} implies that
$f_1$ and $f_2$ can be made disjoint via a homotopy.  Let $N_1=S^1\times
\R \P^9$, let $N_2$ be a $4$--dimensional manifold and let $M^{14}$
be a hyperbolic manifold.  Then one can prove directly that any pair
$(f_1,f_2)$ of maps $f_1\co  N_1\to M^{14}$ and $f_2\co  N_2\to M^{14}$
is homotopic to a pair with disjoint images, and so $\mu(f_1,f_2)=0$. On
the other hand, the Kervaire approach is not applicable if we assume that
$f_1$ is non-trivial on $\pi_1$.  The Hatcher--Quinn~\cite{HQ} Theorem
2.2 is not applicable, since $S^1\times \R \P^9$ is not immersible into
$M^{14}$. (Since $\ov w_6(\R P^9)\ne 0$, neither $\R \P^9$ nor $S^1\times
\R \P^9$ are immersible into the universal cover $\R^{14}$ of $M^{14}$.)
\end{ex}

\begin{proof}
It is clear that if $f_1\in \no$ and $f_2\in \nd$ with $\Im f_1\cap \Im
f_2=\emptyset$ do exist, then $\mu$ is the zero pairing.

Let us show that if $\mu$ is the zero pairing, then such $f_1, f_2$
do exist. Take $\bar f_1\in \no$, $\bar f_2\in \nd$ so that they are
transversal to each other, and hence $\Im \bar f_1$ and $\Im \bar f_2$
intersect at isolated transverse double points. For a double point $d\in
\Im \bar f_1\cap \Im \bar f_2$ we denote by $[d]\in \bor_0(\B)$ the input
of this double point into $\mu([\bar f_1], [\bar f_2])$. The proof is
via induction on the number of double points of $\Im \bar f_1\cap \Im
\bar f_2$.

Take a double point $p\in \Im \bar f_1\cap \Im \bar f_2$. Since $\mu([\bar
f_1], [\bar f_2])=0$, there exists a transversal double point $q\in
\Im \bar f_1\cap \Im \bar f_2$ such that $[p]=-[q]\in \bor _0(\B)$. In
particular we get that the double points $p$ and $q$ correspond to
the same path connected component of $\B$. Put $p_i$ and $q_i$ (for
$i=1,2$)
to be, respectively, the preimages of the double points $p$ and $q$ on
$N_i$ (for $i=1,2$).

Theorem~\ref{pi0} and the definition of $\B$ imply that there exist paths
$\alpha_i\co [0,1]\to N_i$ with $\alpha_i(0)=p_i$ and $\alpha_i(1)=q_i$
for $i=1,2$; and a path $\beta\co [1,2]\to M$ with $\beta(1)=q$,
$\beta(2)=p$ such that the loop $\bar f_i(\alpha_i)\beta$ acts trivially
on the pointed homotopy class of the mapping $(N_i, p_i)\to (M, p)$
for $i=1,2$.  Thus $\bar f_1(\alpha_1)\beta\in \pi_1(M,p)$ commutes with
all the elements of $\Im (\bar f_{1*}\co \pi_1(N_1)\to \pi_1(M))$.

If $\Im \bar f_1 \subset \pi_1(M)$ is an infinite cyclic group or a
trivial group, then, since $\pi_j(M)=0$ for $j=2, \dots,\max(n_1, n_2)+1$,
the elementary obstruction theory implies that $\bar f_1$ is homotopic
to a mapping that passes through a mapping $S^1\to M$. Since $n_1+n_2=m$
and $n_1>2$ we get that after the small perturbation the mapping $S^1\to
M$ does not pass through $\Im \bar f_2$ and thus $\bar f_1$ is homotopic
to a mapping that has disjoint image with $\bar f_2$.

If $\Im \bar f_{1*} \subset \pi_1(M)$ is a group which is not
trivial and is not infinite cyclic, then, since $\pi_1(M)$ is
Preissman (see Definition~\ref{Preissmangroup}), $\bar f_1(\alpha_1)\beta=1\in
\pi_1(M,p)$. Similarly, we get that either $\bar f_2$ is homotopic
to mapping passing through a mapping $S^1\to M$, and hence it is
homotopic to a mapping with the image disjoint from $\Im \bar f_1$;
or $\bar f_2(\alpha_2)\beta$ is also trivial in $\pi_1(M)$. Thus
the paths $\bar f_1(\alpha_1)$ and $\bar f_2(\alpha_2)$ can be
assumed to be homotopic. Thus $\bar f_1(\alpha_1) \bigl( \bar
f_2(\alpha_2)\bigr)^{-1}=1\in \pi_1(M)$ and this loop bounds a disk
$D\subset M$.

For dimension reasons $D$ can be assumed to be embedded and $\Int D$ can be 
assumed to be disjoint from $\Im \bar f_1$ and $\Im \bar f_2$. Since 
$[p]=-[q]\in \bor _0(\B)$, we get that the signs of the intersection points $p$ 
and $q$ of $\Im \bar f_1\cap \Im \bar f_2$ are opposite. Thus we can apply the 
Whitney trick to the points $p$ and $q$ and the disk $D$ and cancel the double 
points. Since $D$ was assumed to be embedded and $\Int D$ was assumed to be 
disjoint from $\Im \bar f_1$ and $\Im \bar f_2$, no new double points of $\Im 
\bar f_1\cap \Im \bar f_2$ will appear during the Whitney trick. 
Thus we managed to decrease the number of double points of $\Im \bar 
f_1\cap \Im \bar f_2$ and we proceed by induction. 
\end{proof}

\end{document}